\input amstex






\font\rm=cmr10 \rm

\font\bf=cmb10
\font\Rm=cmr9 at 11pt
\rm
\font\it=cmsl9 at 10pt
\font\sc=cmr7 
 at 7pt
\def\Sc #1{{\sc \uppercase{#1}}}
\font\Rrm=cmr17 at 16pt
   \font\Rm=cmr12 at 11.5pt

\long\def\Pf{\par\noindent {\it Proof.} }
\def\({\left(}
\def\){\right)}
\def\st{such that }
\def\qed{\hfill$\bullet$\vskip 4pt}

\def\brcs#1{\left\{ #1\right\}}

\def\Log{\text{Log\,}}
\def\iso{\cong}
\def\wrt{with respect to }
\def\:{\,:}

\def\supp{\text{supp}\,}

\def\rk{\text{rank\,}}

\def\ker{\text{ker\,}}

\def\C{\text{\bf C}}

\def\I{\text{I\,}}

\def\R{\text{\bf R}}
\def\N{\text{\bf N}}
\def\Z{\text{\bf Z}}
\def\Q{\text{\bf Q}}

\def\Arrow #1;#2.{#1\:#2 \to }

\def\Set#1#2{\brcs{#1 \left|\vphantom{#1 #2} \right.#2}}



 \def\supp#1{\text{supp\,}#1}
\def\Rrr#1,#2{{\Cal J}_{#1,#2}}
\def\slfrac#1#2{{\raise -.07 ex\hbox{$^{#1}$}}\!/\raise .35 ex \hbox{${}_{#2}$}}
\def\ssf #1/#2{\slfrac {#1}{#2}}

\def\pd #1,#2.{\frac {\partial #1}{\partial #2}}

   \long\def\Lem
#1.#2\par{\vskip4pt{\baselineskip=13pt\font\it=cmsl12 at
11.5pt\Rm
   \noindent {\rm \uppercase{#1}} #2\vskip3pt

   }} 

\long\def\Proclaim #1.#2 \endproclaim{\vskip4pt{\baselineskip=13pt\font\it=cmsl12 at
11.5pt\Rm
   \noindent {\rm \uppercase{#1}} #2\vskip3pt

   }} 

\long\def\remark #1\endremark{\vskip 2pt \noindent {\it Remark\/} #1\par}

\long\def\Sectionhead #1.#2:\par #3{\vskip 4pt \noindent {\bf #1 #2}vskip 2pt\noindent\nospace #3}

\long\def\Title #1\par {\noindent{\Rrm #1}\vskip 9pt}

 \long\def\SubT #1.{\noindent {\it #1\/} } 
 
 \long\def\SecT
#1\par{\vskip 3pt \noindent {\bf #1}\vglue1pt
   \noindent}

\long\def\subtitle #1.{\vskip 2pt \noindent {\it #1}}

\long\def\Rmk#1\par{\vskip 1pt \noindent {\it
Remark.} #1\vskip2pt}

\long\def\Abstract #1\par{{\leftskip= 3 true cm \rightskip = 3 true cm \font\it=cmsl10 \font\rm=cmr10 \baselineskip = 10pt
\parindent=.35 true cm\rm\noindent 
{\it Abstract} #1\vskip 8pt

}}

\long\def\Author #1 \par{\noindent{\it #1}}

\scrollmode\NoBlackBoxes
\magnification=1100

\font\rm=cmr10 \rm

\font\bf=cmb10
\font\Rm=cmr9 at 11pt
\rm
\font\it=cmsl9 at 10pt
\font\sc=cmr7 at 7pt 
 
\def\Sc #1{{\sc \uppercase{#1}}}
\font\Rrm=cmr17 at 16pt
   \font\Rm=cmr12 at 11.5pt

\def\Sc #1{{\sc \uppercase{#1}}}

\long\def\Pf{\par\noindent {\it Proof.} }
\def\({\left(}
\def\){\right)}
\def\st{such that }
\def\qed{\hfill$\bullet$\vskip 4pt}

\def\brcs#1{\left\{ #1\right\}}
\def\Set#1#2{\brcs{#1 \left|\vphantom{#1 #2} \right.#2}}

\def\C{\text{\bf C}}

\def\Log{\text{Log\,}}
\def\wrt{with respect to }
\long\def\Lemma #1. #2\par{\noindent {\Sc  {#1.}} {\Rm
#2}\vskip 2pt}
\def\Arrow #1;#2.{#1\:#2 \to }


\def\R{\text{\bf R}}
\def\N{\text{\bf N}}
\def\Z{\text{\bf Z}}
\def\Q{\text{\bf Q}}
 \def\supp#1{\text{supp\,}#1}
\def\slfrac#1#2{{\raise -.07 ex\hbox{$^{#1}$}}\!/\raise .35 ex \hbox{${}_{#2}$}}
\def\ssf #1/#2{\slfrac {#1}{#2}}


\font\Rrm=cmr17 at 16pt
   \font\Rm=cmr12 at 11.5pt
   \long\def\Lem
#1.#2\par{\write1{#1,
p\,\folio\par}\vskip4pt{\baselineskip=13pt\font\it=cmsl12 \Rm
   \noindent {\rm \uppercase{#1}} #2\vskip3pt

   }} 

   \long\def\Title #1\par {\noindent{\Rrm #1}\vskip 9pt}
   \long\def\SubT #1.{\noindent {\it #1\/} }
   \long\def\SecT #1\par{\vskip 4pt \noindent {\bf #1}\vglue1pt
   \noindent}


\def\oneone{1.1}
\def\onetwo{1.2}
\def\onethr{1.3}
\def\onefou{1.5}
\def\onefiv{1.6}

\def\twoone{2.1}
\def\twotwo{2.2}
\def\twothr{2.3}
\def\twofou{2.4}

\def\throne{3.1}
\def\thrtwo{3.2}
\def\thrthr{3.3}

\def\fouone{4.1}
\def\foutwo{4.2}
\def\fouthr{4.3}
\def\foufou{4.4}
\def\foufiv{4.6}

\def\sixone{6.1}

\def\sevone{7.1}
\def\sevtwo{7.2}
\def\sevthr{7.3}
\def\sevfou{7.4}

\input diagrams

\def\Q{\text{\bf Q}}
\def\Inf{\text{Inf\,}}

\def\paren #1{\/{\rm(}#1\/{\rm)}}

\let\hat=\widehat

\def\Aff{\text{Aff\,}}
\NoBlackBoxes

\def\oneone{1.1}
\def\onetwo{1.2}
\def\onethr{1.3}
\def\onefou{1.4}
\def\onefiv{1.5}

\def\twoone{2.1}
\def\twotwo{2.3}
\def\twothr{2.4}
\def\twofou{2.2}

\def\throne{3.1}
\def\thrtwo{3.2}
\def\thrthr{3.3}

\def\fouone{4.1}
\def\foutwo{4.2}
\def\fouthr{4.3}
\def\foufou{4.4}
\def\foufiv{4.5}

\def\sixone{6.1}
\def\sixtwo{6.2}
\def\sixthr{6.3}

\def\sevone{7.1}
\def\sevtwo{7.2}
\def\sevthr{7.3}
\def\sevfou{7.4}
\def\sevfiv{7.5}

\def\eigone{8.1}
\def\eigtwo{8.2}
\def\eigthr{8.3}
\def\eigfou{8.4}
\def\eigfiv{8.5}
\def\eigsix{8.6}
\def\eigsev{8.7}

\def\I{\text{I\,}}

\let\iso=\cong

\def\tripnorm #1xxx{\left\|\hglue-.2ex\left|#1\right|\hglue-.2ex\right\|}

\Title Non-direct limits of simple dimension groups with finitely many pure traces

\Abstract There exist simple dimension groups which cannot be expressed as a direct limit of simple, or even approximately divisible  dimension groups, each with finitely many pure traces, and we can specify its infinite-dimensional Choquet simplex of traces; a more drastic property is noted. On the other hand, a very easy argument shows that if $G$ is a $p$-divisible simple dimension group (for some integer $p>1$), then it can be  expressed as such a direct limit. We also enlarge the class of initial objects for AF (and slightly more general) C*-algebras.

\noindent  {\it David Handelman}%
\plainfootnote{$^1$}{\hglue -.5em Supported in part by a Discovery grant from NSERC.}

\noindent Thinking about properties of traces on dimension groups (see for example, [BeH]), especially simple dimension groups, I realized that it would have been nice to be able to reduce to simple dimension groups with finite pure trace space, or better, to simple dimension groups of finite rank (these automatically have finitely many pure traces), or better still, to simple dimension groups  finitely generated as abelian groups.

\noindent This suggests three conjectures:

\Lem Conjectures. Every noncyclic simple dimension group is a direct limit (that is, over a directed set, and with positive maps) of simple dimension groups
\itemitem{(a)} that are free of finite rank;
\itemitem{(b)} that have finite rank;
\itemitem{(c)} that have only finitely many pure traces.

\noindent Obvously, the truth of conjecture (a) would imply that of (b), which would in turn imply that of  (c). In fact, it turns out that (c) is false (so that (a) and (b) are too), via a reasonably well-known example, and not only that, it is false in an extreme way: there is a noncyclic simple dimension group $R$ (in fact, an ordered ring) for which there are no nonzero positive group homomorphisms  $\Arrow \alpha; G.R$ for any approximately divisible dimension group $G$  with only finitely many pure traces. (All noncyclic simple dimension groups are approximately divisible, so this is a drastic way of showing $R$ cannot be a direct limit of those with only finitely many pure traces.)

This contrasts with the corresponding (dual) question for Choquet simplices: a fundamental  result of Lazar \& Lindenstrass [LL] is that every metrizable Choquet simplex is an inverse limit (over the positive integers) of finite-dimensional simplices. As pointed out in [G], a quick proof of this is derivable using dimension groups, from the main result of [EHS]. 

However, there is an easy result that yields many examples that are such direct limits, even for  case (b). Case (a) is somewhat problematic, but there is an example later, $\Z[\sqrt 2][x]$ equipped the strict ordering as functions on $[\slfrac 13, \slfrac 23]$, which is a direct limit for case (a) (the pure trace space is a continuum, so the result is not entirely trivial).

All groups (and partially ordered groups) appearing here are torsion-free abelian; {\it free\/} means free as an abelian group. Recall that the {\it rank\/}\plainfootnote{$^2$}{\hglue -.5em Some authors have very unfortunately defined the {\it rank of a dimension group\/} to be the width of the minimal Bratteli diagram realizing it. This is different from the rank of the underlying group, and should be called the minimal width.}
 of a (torsion-free) abelian group $J$ is the dimension over $\Q$ of the rational vector space $J \otimes_{\Z}\Q$. All partially ordered groups appearing here are unperforated.

The first section gives fairly easy examples of simple dimension groups that can be expressed as a direct limit of simple dimension groups with finitely many pure traces; this leads to a definition of pro-finite dimensional (for a somewhat larger collection of dimension groups with order unit) and the class consisting of them is denoted $\Cal H$. The second section is devoted to construction of a single example of what we call (strongly) anti-finite dimensional (Anti-FD), a simple dimension group $R$ \st for all $G \in \Cal H$, the only positive group homomorphism $\Arrow \phi;G.R$ is zero ($R$ has an even stronger property than this). 

The third section shows that Anti-FD simple dimension groups can be constructed to have any Choquet simplex as their trace space (and if the Choquet simplex is metrizable, the dimension group can be chosen to be countable as well). Finally, the fourth section provides a characterization of Anti-FD; it boils down to triviality of the infinitesimal subgroup and an elementary topological constraint, that every finite rank subgroup be discrete \wrt the norm obtained from the affine representation. 

Section 6 deal with examples arising from actions of tori on UHF algebras. Sections 7 and 8 show that a small class of the latter are initial objects in the appropriate category, enlarging on   work of Elliott and R\o rdam [ER]. 

\SecT 1 Pro-finite dimensional dimension groups
 
The following formalizes a well-known construction.
 Let $G$ be a partially ordered abelian group (not necessarily a dimension
group) having an order unit $u$. We define the {\it simplification of
$G$}, $G^s$, to be the same underlying group $G$, but with positive cone
consisting of the set of order units of $G$ (denoted $G^{++}$) together
with $0$. This puts a generally coarser ordering on the group, although
the identity map $G^s \to G$ induces a natural affine homeomorphism
$S(G,u) \to S(G^s,u) $ (since $G^+ + G^{++} = G^{++}$, and the order units
of $G$ and of $G^s$ are identical). The simplification of an ordered group
is simple, and obviously, if $G$ is already simple, it equals its
simplification. Although not formally named, this process occurs in many (well, some) papers on dimension groups. 
 
In general, if $G$ is a dimension group, $G^s$ need not be: for any $n >1$, take $G= \Z^n$
with the usual ordering, so that $G^s$ is $\Z^n$ with the strict ordering,
that is, its positive cone consists of $0$ together with the set of
$n$-tuples, $(k(1), \dots, k(n))$
\st $k(i) > 0$ for all $i$. This admits discrete traces, so cannot be a simple dimension group [G]. 
 
A dimension group $G$ with order unit is  {\it approximately divisible\/}
if for all pure traces $\tau$, $\tau(G)$ is a dense subgroup of $\R$. This
is compatible with other definitions of approximately divisible because of
the following result.
 
\Lem Theorem \oneone. [GH] A dimension group $G$ with order unit is
approximately divisible iff its image in $\Aff S(G,u)$ is norm dense for
one (hence for all) order unit(s) $u$. All noncyclic simple dimension
groups are approximately divisible.
 
In particular, if $A$ is a unital AF C*-algebra, then K${}_0 (A)$ is
approximately divisible if and only if $A$ has no finite-dimensional
representations. This combined with the fact that a dense subgroup of
$\Aff K$, equipped with the strict ordering, is a dimension group iff $K$
is a (Choquet) simplex, yields the following.
 
\Lem Corollary \onetwo. Let $(G,u)$ be a partially ordered abelian group with more
than one pure trace. Then its simplification is a dimension group if and
only if (a) $S(G,u)$ is a Choquet simplex and (b) the image of $G$ is norm
dense in $\Aff S(G,u)$. In particular, if $G$ is an approximately
divisible dimension group, then $G^s$ is a dimension group.
 
Now suppose the partially ordered group $G$ with order unit $u$ can be
expressed as a direct limit of unperforated partially ordered groups, in
particular, $G = \lim \Arrow \phi_n ;G_n.G_{n+1}$ (the following remarks
also apply when uncountable directed sets are used, however, not only does the notation
become tedious, but uncountable direct limits are not very important in this context), where $\phi_n$ are order preserving (that is, {\it
positive\/} group homomorphisms). Suppose each $G_n$ has an order unit
$u_n$ \st $\phi_n (u_n)$ is an order unit of $G_{n+1}$.
 
Then $\phi_n (G_n^{++}) \subseteq G_{n+1}^{++}$. To see this, note that $x
\in G_n^{++}$ iff there exists an integer $k$ \st $u_n \leq k x$; hence
$k\phi_n(x) \geq \phi_n (u_n)$, so that $k\phi_n(x)$ is an order unit in
$G_{n+1}$, and by unperforation, $\phi_n(x)$ is thus an order unit. Hence,
the $\phi_n$ are positive homomorphisms $\Arrow \phi_n; G_n^s. G_{n+1}^s$
between their simplifications.
 
The identity mapping sends $\lim \Arrow \phi_n; G_n^s. G_{n+1}^s$ to $G =
\lim  \Arrow \phi_n; G_n. G_{n+1}$, and is obviously order preserving
(since the positive cones in the simplified groups are contained in the
positive cones of the originals). If we now assume that $G$ itself is
simple, then this map is an order-isomorphism: if $g$ is in $G^+\setminus
\brcs{0}$, then $g \in G^{++}$, so there exists a positive integer $k$ \st
$k g \geq u$; hence there exists $n$ \st $g$ is represented by $[h,n]$ and
$u$ is represented by $[v,n]$ with $v,h \in G_n^+$ and $v \in G_n^{++}$,
and $k h \geq v$ in $G_n$. Hence $h$ is an order unit of $G_n$ and $[h,n]$
represents $g$. Thus $h$
 is positive as an element of $G_n^s$, so everything that is positive in
$G$ is represented as the image of a positive element in $\lim \Arrow
\phi_n; G_n^s. G_{n+1}^s$.
We have thus proved,
 
\Lem Lemma \onethr. Suppose that $G = \lim \Arrow \phi_n; G_n. G_{n+1}$ is a
simple partially ordered group, each $G_n$ is unperforated, and each $G_n$
contains an order unit $u_n$ \st $\phi(u_n)$ is an order unit in
$G_{n+1}^s$. Then $G$ is also the ordered direct limit of simple
unperforated groups, $\lim \Arrow \phi_n; G_n^s. G_{n+1}^s$.
 
If we begin with a not necessarily simple ordered group $G = \lim \Arrow
\phi_n; G_n. G_{n+1}$ but with the same hypotheses on $G_n$ and $\phi_n$,
then we quickly see that $G^s = \lim \Arrow \phi_n; G_n^s. G_{n+1}^s$ by
essentially the same argument.
 
If $G$ is a dimension group, we can write $G= \lim \Arrow A_n; \Z^{k(n)}.
\Z^{k(n+1)}$. To satisfy the order unit condition, we telescope and delete
rows if necessary so that the matrices $A_n$ have no zero rows. Then $G^s
= \lim \Arrow A_n; (\Z^{k(n)})^s.( \Z^{k(n+1)} )^s$. This is not very
illuminating, since $(\Z^{k(n)})^s$ is not a dimension group (unless $k(n)
=1$). However, if $G$ is a $\Z[1/p]$-module for some integer $p > 1$ (that
is, $G$ is $p$-divisible), then $G \iso G \otimes \Z[1/p]$ (ordered tensor
product over $\Z$), and thus we can write
$G \iso  \lim \Arrow A_n; \Z[1/p]^{k(n)}. \Z[1/p]^{k(n+1)}$. Now
$\Z[1/p]^k$ is dense in $\R^k$, and so its simplification is a dimension
group, obviously simple.

A dimension group $G$ is {\it pro-finite dimensional\/} (pro-fd) if it can be realized as a direct limit (over a directed set) of positive group homomorphisms $\Arrow \phi_{\alpha, \beta}; G_{\alpha}. G_{\beta}$ where each $G_{\alpha}$ is an approximately divisible dimension group with order unit, $u_{\alpha} $, \st $\phi_{\alpha,\beta}(u_{\alpha})$ is an order unit in $G_{\beta}$, and each $\partial_e S(G_{\alpha},u_{\alpha})$ is finite (finite-dimensional trace space). Necessarily, a pro-FD dimension group has an order unit. Most of the time, specifically when $G$ is countable, the index set can be taken to be the set of positive integers, in which case the notation simplifies to $\Arrow \phi_i; G_i. G_{i+1}$. The class of pro-FD dimension groups will be denoted ${\Cal H}$. Rather surprisingly, not all simple noncyclic dimension groups belong to $\Cal H$.

We have seen that if $G$ is simple and in $\Cal H$, then $G$ can be represented as a limit of simple dimension groups each with finite-dimensional trace space.
 
\Lem Corollary \onefou. If $G$ is a dimension group that is $p$-divisible for some
integer $p > 1$, then $G^s$ is a direct limit of simple dimension groups each 
with finitely many pure traces. In particular, every simple  $p$-divisible
dimension group is a direct limit of simple dimension groups, each having
only finitely many pure traces, and is thus pro-fd.
 
A $p$-divisible dimension group is automatically approximately divisible.
A similar construction occurs if $G$ can be represented as a limit of the
form $R^{k(n)}\to R^{k(n+1)}$ where $R$ is a noncyclic ordered subring of
the reals (for example, $R = \Z[\sqrt 2]$. If $G$ is simple, $G$ is again a direct
limit of simple dimension groups each with only finitely many pure traces.

\Lem Corollary \onefiv. If $G$ is a simple dimension group that is $p$-divisible for some prime $p$, then $G$ is a limit of simple dimension groups of finite rank.

\comment
\Pf (Embarassingly easy.) Write $G $ as a limit of simplicial groups in the usual way; because $G$ is simple, we may telescope, and so assume the matrices $\Arrow A_n;\Z^{k(n+1)}.\Z^{k(n)}$ are strictly positive. Now assuming $G$ is $p$-divisible, then $G \iso G\otimes_{\Z} \Z[1/p]$ not only as groups, but as partially ordered abelian groups (this is trivial, and follows easily from unperforation).

Set $H_n = \Z[1/p]^{k(n+1)} = \Z^{k(n+1)} \otimes \Z[1/p]$. If we impose the usual direct sum ordering on $H_n$, then the limit (using the same matrices) is $G$ again; if instead, we impose the strict ordering on $H_n$, then each $H_n$ is a simple dimension group, and since the matrices $A_n$ are strictly positive, they implement positive $\Z[1/p]$-homomorphisms, and so we can take the ordered direct limit. The underlying group is the same, but the new positive cone is contained in the original (using the ordinary ordering); however, if some nonzero element is nonnegative in the original ordering, then as $G$ is simple, it is represented at some level by a strictly positive column. Hence the orderings are the same, and we have realized $G$ as a direct limit of simple dimension groups, each of finite rank.
\qed

\endcomment

This likely extends to groups satisfying the following property. Say an element $g$ of a torsion-free abelian group is {\it i-divisible\/}  ({\it i} for {\it infinitely\/}) if there exist infinitely many positive integers $n$ \st the equation   $g = nx_n$ can be solved (with $x_n$ in $G$). Now say the abelian group $G$ is {\it wi-divisible\/} if every element is a sum of i-divisible elements, that is, the i-divisible elements span the group. It is plausible that $p$-divisible can be weakened to wi-divisible in the statement of the result.

As an example, let $H_n = \oplus_{i=1}^{k(n)} U_{i,n}$ where each $U_{i,n}$ is a noncyclic subgroup of $\Q$, and let $(A_n)$ be a sequence of strictly positive matrices, with $A_n$ being of size $k(n+1) \times k(n)$ and having integer entries. We can take $G = \lim \Arrow A_n;H_n.H_{n+1}$. The outcome is a simple dimension group that is $wi$-divisible; uninterestingly, because the $A_i$ are strictly positive, there is a noncyclic subgroup $U$ of $\Q$ for which we can rewrite the terms as $U^{k(n)}$, so that the  group is i-divisible.  If we impose the strict, rather than the ordinary direct sum ordering on each $H_n$, then each $H_n$ is now a simple dimension group, and the $A_n$ are still positive homomorphisms between them (this is why we insisted the $A_n$ have no zero rows; if there were any zero rows, some positive elements would be sent to a non-positive elements), and as in the argument above, simplicity of $G$ entails that the limit is just $G$ with its original ordering. Hence $G$ is a limit of simple dimension groups of finite rank, specifically the $H_n$ with the strict ordering.

On the other hand, $G = \Z[\slfrac12] \oplus \Z[\slfrac13]$ with the strict ordering (as a subgroup of $\R^2$) is wi-divisible without being i-divisible. 

\SecT 2 Drastic example

At the opposite extreme to pro-fd, are two properties of dimension groups defined below. An approximately divisible dimension group $(G,u)$ with order unit is {\it anti-finite dimensional\/} (anti-fd) if there are no nonzero positive group homomorphisms $\Arrow \phi; H.G$ for any approximately divisible dimension group $H$ having  finite-dimensional pure trace space (equivalently, for any  $ H \in \Cal H$). This is a much stronger property than merely not being a limit of simple dimension groups with finitely many pure traces. That simple  anti-fd dimension groups even exist is somewhat surprising; however, it turns out that for every Choquet simplex $K$, there exists one with trace space $K$ (and if $K$ is metrizable, we can find a countable one).
Actually we show there exist lots of approximately divisible dimension groups with an even stronger property.

Let $(G,u)$ be a partially order abelian group with order unit, and form the representation  $G \to \Aff S(G,u)$, given by $g \mapsto \hat g$, where $\hat g (\tau) = \tau(g)$ as usual. This induces a pseudo-norm on $G$ arising from the supremum norm on $\Aff S(G,u)$; it can also be characterized purely in terms of the ordering on $G$. If the ordered group is unperforated, then the kernel of the representation consists of the infinitesimals, $\Inf G$. We say a group homomorphism $\Arrow \phi;G.H$ between partially ordered groups with order unit (but not necessarily sending order units to order units) is {\it continuous\/} if it is continuous \wrt the pseudo-norms on $G$ and $H$; equivalently, $\phi(\Inf G) \subset \Inf H$ and the induced map $\Arrow \overline{\phi}; G/\Inf G. H/\Inf H$ is continuous \wrt to the norms on the groups induced by their affine representations.

When $G$ has finite-dimensional trace space, every continuous group homomorphism $G \to H$ is bounded, hence will extend to a norm-continuous map between their completions. A continuous (or boundedÑsee [G, section 7] for a discussion of bounded group homomorphisms on dimension groups) group homomorphism from a dense subgroup of $\R^n$ to any Banach space automatically extends to a bounded linear function from $\R^n$ to the Banach space (since continuous maps need not send Cauchy sequences to Cauchy sequences, this is not completely obvious; however, if $x_n \to x/k$ where $x_n$ and $x$ are in the dense subgroup (and $k$ is an integer), then $kx_n \to x$, so continuity implies $\alpha(kx_n) \to \alpha(x)$, and this  implies extension to a rational subspace, to which the standard method of showing continuity implies boundedness for linear transformations applies.) 
Every positive group homomorphism is automatically continuous, so the following even more  drastic property implies anti-finite dimensionality.

We say a dimension group with order unit $R$ is {\it strongly anti-finite dimensional\/} (Anti-FD; we use three upper case letters to distinguish it from anti-fd and  Anti-fd, the latter being anti-fd at the beginning of a sentence) if for every $H$ in $\Cal H$, every continuous group homomorphism $\Arrow \phi; H.R$ is zero. This condition forces $\Inf R = \brcs{0}$ for trivial reasons. A few examples of simple anti-fd dimension groups with no infinitesimals that are not Anti-FD are given in section 4.

While characterization of anti-fd dimension groups is a little complicated, the characterization of Anti-FD dimension groups (within the class of approximately divisible dimension groups with order unit) is easy: $\Inf R = \brcs{0}$ and every finite rank subgroup is discrete (in the norm topology), Corollary \foufiv. In this section we show sufficiency, then give more examples (having as trace space, all Choquet simplices) in the next, and finally show necessity in the last section. The most interesting feature is that Anti-FD dimension groups exist. This means that $\Cal H$ is far from exhaustive, even for simple dimension groups.

Now we prepare to provide an example of a simple dimension group that cannot be obtained as a direct limit of simple dimension groups, each with finitely many pure traces, in fact is Anti-FD, and in particular, is a counter-example to (c).

For a pseudo-normed abelian group $H$, let $\psi$ denote the map to its completion; that is, first factor out the elements of norm zero, and complete the resulting normed abelian group. Of course, the completion need not be a vector space (e.g., if the norm is discrete).

\Lem Proposition \twoone. Suppose that $H$ is a pseudo-normed group \st the completion of $\psi(H)$  a finite dimensional real vector space, and $R$ is a normed abelian group with the following properties:\hfill\break
{(0)} the completion of $R$ is a Banach space \hfill\break
{(a)} every countable subgroup of $R$ is free as an abelian group \hfill\break
{(b)} every finitely generated subgroup of $R$ is discrete.\hfill\break
\noindent Then every  continuous group homomorphism $\Arrow \phi; H.R$ is zero.

\Pf We immediately reduce to the case that $H$ is a dense subgroup of a finite dimensional real vector space (so we can dispense with $\psi$). Now consider the subgroup $\phi(H) $ of $R$. If $\phi(H)$ is nonzero and has finite rank (as an abelian group), it is countable, as a subgroup of a free group, it is itself free; a free group of finite rank is finitely generated. Hence $\phi(H)$ is discrete. There thus exists a continuous linear functional on the completion of $R$, $\rho$, \st $\rho(\phi(H)) = \Z$.

Then $\Arrow\upsilon:= \rho \circ \phi;H.\R$ is a continuous group homomorphism. Since the completion of $H$Ñactually pseudo-completion, but $\psi$ kills the subgroup consisting of those elements of $H$ with zero norm, so we might as well assume $H$ is dense in $\R^n$Ñis  a finite dimensional vector space, $\upsilon$ extends to a bounded linear functional on the completion. But its restriction to (the image of) $H$ has discrete range, which is impossible since $H$ is dense in a vector space.

We are thus reduced to the case that $\phi(H)$ has infinite rank. Suppose the (real) vector space dimension of the completion of $H$ is $n$. Every abelian group of infinite rank contains a subgroup of every finite rank. Hence there exists a rank $n+1$ subgroup $J$ of $\phi(H)$. As a subgroup of $R$, it is free, and thus free on $n+1$ generators, call them $\phi(h_i)$. By (c), the group they generate is discrete, and thus there exist continuous linear functionals $\alpha_i$ on $B$ (and thus continuous real-valued group homomorphisms on $R$) \st $\alpha_i (\phi(h_j)) = \delta_{ij}$ (Kronecker delta).

Then $\brcs{\upsilon_i:= \alpha_i \circ \phi}$ is a collection of continuous additive real-valued group homomorphisms on $H$, each of which extends to a bounded linear functional on the completion. However, $\upsilon_i(h_j) = \delta_{ij}$ entails that the
set of bounded linear functionals is linearly independent (over the reals of course), so that the dimension of the dual space of the completion, and thus of the completion itself, exceeds $n$, a contradiction.
\qed

Properties (a) and (b) together can be restated as a single property (of normed abelian groups),
\item{(¥)} every finite rank subgroup of $R$ is discrete.

\noindent To see the equivalence, note that a finite rank discrete group is free [actually, a result due to Steprans [S], generalizing a result of Lawrence [L], asserts that {\it every\/} discrete normed abelian group is free]; by [Gr, Theorem 137, p\,101], a countable group for which every finite rank subgroup is free, is itself free. And of course, if $R$ is a free abelian group, then every subgroup is free.

It is convenient, however, to separate (a) and (b), in view of the closely related examples we obtain. 

If $\brcs{x_i}_{i=1}^n$ is a finite and (real) linearly independent subset of a Banach space, then the abelian group the set generates, $\sum x_i \Z = \oplus x_i \Z$ not only is free (as an abelian group) on $\brcs{x_i}$ but is discrete: the map $x_i \mapsto e_i$ (where $e_i$ run over the standard basis of $\R^n$) extends to an automatically continuous and real linear isomorphism (with continuous inverse) $\sum x_i \R \to \R^n$ that sends $\sum x_i \Z$ to the usual copy of $\Z^n \subset \R^n$.

\Lem Example \twofou. A simple Anti-FD dimension group.

\noindent Set $R = \Z[x]$, the ring of polynomials with integer coefficients. Let $I$ be a closed real interval (of nonzero length) that contains no integers. Then a 1925 result due to Chlodovsky [C] (as cited in [F]), says that $R$ is dense in $C(I,\R)$ \wrt the supremum norm on $I$. In particular, if we impose the strict ordering on $R$ (that is, $f \in R^+\setminus{0}$ iff $f$ is strictly positive as a function on $I$), then $R$ is a  simple dimension group.   We note that $R$ is an ordered ring with $1$ as order unit, and its pure trace space consists of the point evaluations at members of $I$, so its pure trace space is an interval.

Moreover, $R = \sum x^i \Z$ expresses $R$ as  free  on $\brcs{x^i}$. Finally, every finitely generated subgroup of $R$ is contained in $L= \sum_{i=0}^n x^i \Z$ for some $n$, so to prove the subgroup is discrete, it suffices to show $L$ is discrete. However, $\brcs{x^i}$ is real linearly independent  as a subset of $C(I,\R)$ (a polynomial with real coefficients that vanishes on $I$ is zero), and thus $L$ is discrete. So $R$ is a simple dimension group satisfying conditions (0, a, b) of Proposition \twoone, and is thus Anti-FD. 

The affine representation of $R $ (using $1$ as an order unit) is simply the inclusion $R \subset C(I,\R)$ and its pure trace space  is naturally homeomorphic to $I$). \qed

If $H$ is an approximately divisible dimension group with order unit, then the image of $H$ is dense in its affine representation. Hence if $H$ is approximately divisible and has only finitely many pure traces, then the natural pseudo-norm from the affine representation has Euclidean space as its completion. The kernel of the affine repesentation consists of infinitesimals, which is exactly the subset to be factored out in constructing the completion.
Any positive homomorphism between partially ordered abelian groups with order unit is automatically continuous \wrt the pseudo-norm topologies on each.

\comment
Let $\Cal H$ denote the class of direct limits of approximately divisible
dimension groups with only finitely many pure traces, \st the connecting maps send order units to order units. Then the members of $\Cal H$ are themselves
approximately divisible dimension groups with order unit, and as such
acquire a pseudo-metric topology from the map $H \to \Aff S(H,v)$ (for any
order unit $v$ of $H$).
\endcomment

An abelian group homomorphism $\Arrow \phi; H.J$ where $H$ is a pseudo-normed abelian group and $J$ is a normed abelian group is {\it weakly continuous\/} if for every continuous group homomorphism $\Arrow \rho;J.\R$, the composition $\Arrow \rho\circ \phi; H.\R$ is continuous. The proof of Proposition \twoone\ only requires that $\phi$ be weakly continuous. However, when the pseudo-completion of  $H$ is a finite-dimensional real vector space, weak continuity implies continuity anyway. In the following result, dealing with maps from members of $\Cal H$, we have two ways of proceeding: directly, using weak continuity of the restriction to the image of something with finite dimensional completion, or using that the restriction of a weakly continuous group homomorphism to one of the constituents in the direct limit is automatically continuous.
 
With $R$ the example of \twofou\ (or any other anti-fd simple dimension group), we have the following. Of course, $R$ itself
has the metric topology from the sup norm on the interval, and this
coincides with the metric (not just pseudo-metric) obtained from its
affine representation (which would be as a dense subring of $\Aff \Cal
M^+(I) = C(I,\R)$, $\Cal M^+$ being the collection of probability
measures on $I$). 
 
\Lem Lemma \twotwo. Let $H \in \Cal H$. If $\Arrow \phi;H. R$ is a weakly continuous group
homomorphism, then $\phi =0$. This applies automatically
if $\phi$ is an order-preserving group homomorphism.

\comment
We say a dimension group with order unit $(G,u)$ has the {\it anti-finite dimensional property\/} if for every approximately divisible dimension group $(H,v)$ with order unit \st $S(H,v)$ is finite dimensional, there are no  positive homomorphisms $\Arrow \phi; H.G$ \st $\phi(u)$ is an order unit of $H$. We say $(G,u) $ has the {\it strong anti-finite dimensional property\/} if there are no continuous group homomorphisms $\Arrow \phi; H.G$ for the same class of $H$s. Proposition \twoone\ gives sufficient conditions for a dimension group to be strongly anti-finite dimensional, in particular $R$ has the strong anti-finite dimensional property. (Whether there exist simple dimension groups with the anti-finite dimensional property, but not the strong anti-finite dimensional property, is open.)

\Lem Lemma. If $G$ is an approximately divisible dimension group of finite rank, then there is no nonzero positive group homomorphism $\Arrow \phi; G.R$.

\Pf Since $\phi(G)$ is a finite rank subgroup of  $R$, a free abelian group on $\brcs{x^i}_{i\geq 0}$, there exists $n$ \st $\phi (G) \subset \sum_{i=0}^n x^i \Z$. Since $\brcs{x^i}_{i=0}^n$ is linearly independent as a subset of $C([1/2,3/4],\R)$, the latter is discrete in the sup norm. Hence there exists a bounded linear functional on $C =C([1/2,3/4],\R)$, $\alpha$ \st $\alpha(\sum_{i=0}^n x^i \Z) =\Z$. We can write $\alpha = s -t$ where $t$ are positive (unnormalized) linear functionals on $C$. The composite maps $\sigma = s \circ \phi$ and $\tau = t \circ \phi$ are thus traces on $G$, and $(\sigma-\tau)(G) = \Z$. However, being approximately divisible, $G$ has dense range in its affine function representation, in particular, its range under any difference of unnormalized traces must be dense in $\R$, a contradiction.
\qed

In particular, a noncyclic simple dimension group is approximately divisible, so we obtain the a counter-example to (b). But the result is better, since it applies to approximately divisible dimension groups (but that turns out to be a red herring).
\endcomment

This contrasts with $\Z[\slfrac12][x]$ and $\Z[\sqrt 2][x]$, each equipped with the strict ordering from restriction to the interval;  $\Z[\slfrac12][x]$ is 2-divisible, so is even a direct limit of simple dimension groups of finite rank. And $\Z[\sqrt 2][x]$ is order isomorphic to the ordered tensor product $R \otimes_{\Z} \Z[\sqrt 2]$ (where $\Z[\sqrt 2]$ is given the total ordering inherited from $\R$: we can write $R = \lim \Arrow A_n ; \Z^{k(n)}. \Z^{k(n+1)}$ where the $A_n$ are strictly positive matrices (since $R$ is simple, this can be arranged by telescoping), and since $R$ and thus $R \otimes \Z[\sqrt 2]$ is simple, the latter is the limit $\Arrow A_n \otimes 1; \(\Z[\sqrt 2]^{k(n)}\)^s. \(\Z[\sqrt 2]^{k(n+1)}\)^s$, each direct sum equipped with the strict ordering. Every $\(\Z[\sqrt 2]^{k(n)}\)^s$ is a simple dimension group (from being dense in $\R^{k(n)}$). In this latter example, the underlying group is free, and it is a direct limit of simple dimension groups that are finitely generated and free. In particular, both $\Z[\slfrac12][x]$ and $\Z[\sqrt 2][x]$ equipped with the strict ordering from some closed interval in $(0,1)$ are simple dimension groups in $\Cal H$, but their intersection, $R =\Z[x]$ with the strict ordering on the same interval, is a simple dimension group which has no nonzero incoming weakly continuous group homomorphisms from any element of $\Cal H$. 

The ring $\Z[\sqrt 2][x]$ with the strict ordering satisfies (a)
 of Proposition \twoone, but not (b); on the other hand, $\Z[\slfrac 12][x]$ satisfies (b) but not (a). Both satisfy (0).

We actually have uncountably many choices for $I$ that give rise to nonisomorphic simple dimension groups: for example, the unordered pair consisting of left and right endpoints of the interval, $I = [a,b]$, is topologically determined from the pure trace space (since they are the only two points with no neighbourhoods homeomorphic to an open interval), and the corresponding value groups, $\brcs{\Z[a], \Z[b]}$, viewed as subgroups of the reals, is an order-theoretic invariant of $\Z[x]$ with the strict ordering coming from $I$ (and there are uncountably many different order isomorphism classes of $\Z[a]$). In the next section, we will show that arbitrary Choquet simplices can be realized as the trace space of  simple Anti-FD dimension groups.

The argument in Proposition \twoone\ was divided in two parts; first, dealing with subgroups of finite rank, then with those of infinite rank. Had the following  notion (not   good enough to be a conjecture) been true, the second argument would have been unnecessary.\vskip4pt

\itemitem{({\it Notion\/})} If $G$ is a dense subgroup of $\R^n$, there exists a finite rank subgroup of $G$ that is also dense.\vskip 4pt

\noindent This is plausible but false; this phenomenon accounts for problems in  characterizing anti-fd dimension groups, as discussed in section 4.

\Lem Example \twothr. A countable dense subgroup of $\R^2$ which contains no dense subgroups of finite rank. This yields a simple dimension group with exactly two pure traces to which there are no positive homomorphisms from any noncyclic simple dimension group of finite rank.

\noindent Let $\brcs{\alpha_0 = 1, \alpha_1, \alpha_2, \dots}$ be an infinite set of real numbers that is linearly independent over the rationals. For $i = 0,1,2,\dots$ set $x_i = (\alpha_i, 2^{-i}) \in \R^2$. Let $G = \sum x_i \Z \subset \R^2$. By examining the first coordinates, we see that $G$ is a free abelian group on $\brcs{x_i}$, obviously of infinite rank.

We observe that $2x_1 - x_0 = (2\alpha_1 -1,0) $ and $4x_2 - x_0 = (4\alpha_2 -1,0)$. Since $\brcs{2\alpha_1 -1, 4\alpha_2 -1}$ is linearly independent over the rationals, we see that $\R \oplus 0$ is contained in the closure of $G$. Thus for all $i$, $(0, 2^{-i})$ is contained in the closure of $G$, so that $0 \oplus \R$ is also contained in the closure, and thus, as the closure is a group, $\R^2$ is the closure. So $G$ is dense in $\R^2$.

Now suppose $H$ is a subgroup of finite rank. Since $G$ is free, $H$ must be free, and as its rank is finite, it is finitely generated as an abelian group. Hence there exists $n$ \st $H \subseteq \sum_{i=1}^n x_i \Z$. The restriction of the second coordinate functional has values in $2^{-n}\Z$, which is discrete. In particular, $H$ has a trace with discrete range, so cannot be the image of a noncyclic simple dimension group (even dropping the requirement that the image be dense in $G$). 

Impose the strict ordering on $G$; it is the desired simple dimension group; all noncyclic simple dimension groups have dense range in their affine representation.

As a special case, let $\alpha$ be a real transcendental  number, and set $\alpha_i = \alpha^i$. The resulting $G$ is simply the ring $\Z[x]$ equipped with the ordering derived from  the two ring homomorphisms  $x \mapsto \alpha$ and $x \mapsto \slfrac 12$; the image is dense in $\R^2$, making $\Z[x]$ into a simple dimension group (and an ordered ring), with these two homomorphisms as the only pure traces (up to normalization at $1$). Since $\alpha$ is transcendental, the map $x \mapsto \alpha$ has zero kernel. This example is simply $\Z[x]$ with a much larger positive cone than that of $R$ in Example \twofou. 
\qed

\SecT 3 Other Choquet simplices

In the construction of the  Anti-FD simple dimension group $R=\Z[x]$ with the strict ordering from restriction to $I$, the set $I$ is a closed interval; everything works just as well if we let  $I$ be any infinite compact subset of the open unit interval (since a polynomial is determined by its values on any infinite set). As we can embed a Cantor set in the open interval (for example, truncate the usual Cantor set at both ends), we can even arrange that the pure trace space of $\Z [x]$ equipped with the strict ordering be a Cantor set, hence totally disconnected; the result will still be a simple dimension group satisfying (0), (a), and (b) of Proposition \twoone. However, we can go much farther. 

We show that for every (metrizable) Choquet simplex $K$, there exists a (countable) simple dimension group with order unit $(G,u)$ whose trace space is $K$ that satisfies conditions (a) and (b) of the proposition, and therefore has the strong anti-finite dimensional property.

\comment
\Lem Lemma. Let $X$ be a compact Hausdorff space. Then $C(X,\R)$ contains
a linearly independent set $\brcs{f_{\alpha}}$ \st $\sum f_{\alpha }\Z$ is
a dense subgroup and a ring.
 
\Pf Every compact set is embeddable in a product of copies of the unit
interval. Hence we may find an embedding in $Y:= \prod_{\beta \in \Omega}
[
\slfrac13,\slfrac23]_{\beta}$. Let $\brcs{x_{\beta}}$ be correspondingly
indexed variables, and set $\Z[\Omega]:= \Z [x_{\beta}]_{\beta \in
\Omega}$, the integer polynomial ring variables indexed by $\Omega$. We
can rewrite this as the tensor product $\otimes_{\Omega} \Z[x_{\beta}]$,
each one of which has dense image in
$C([\slfrac13,\slfrac23]_{\beta},\R)$, so that on viewing $\Z[\Omega]$ as
functions on $Y$, we have that it has dense image. Now Tietze's extension
theorem yields that the restriction of $\Z[\Omega]$ to the image of $X$.
Hmmm, not quite only works if the restriction to $X$ is one to one... on
polynomials.
\endcomment
 
The following is a minor variation on well-known result.
 
\Lem Lemma \throne. Let $B$ be a separable  infinite-dimensional real  Banach
space, and let $Q$ be a countable linearly independent subset. Then there exists a countable dense set $P$ \st $P \cap Q = \emptyset$ and $P \cup Q$ is linearly independent.
 
\Pf Since $B$ is separable, there exists a  countable dense set $\brcs{p_i}_{i=1}^{\infty}$. Since $B$ is complete and
not finite-dimensional, its dimension (as a real vector space) is
uncountable. Select a basis for $B$ (as a real vector space), $A =
\brcs{e_{\alpha}}$. Every element of $B$ is uniquely representable in the
form $b = \sum \lambda_{\alpha} e_{\alpha}$ with $\lambda_{\alpha} \in \R$
where all but finitely many $\lambda_{\alpha} $ are zero; define the {\it support of $b$,} given by $\supp b
= \Set{e_{\alpha} \in A}{\lambda_{\alpha} \neq 0}$.
 
Let $J_0$ be the union of the supports of all the elements of $Q$; this is countable. Let $J_1$ be the union of the supports of all the elements of $\brcs{p_i}$. Since $J:= J_1 \cup J_0$ is countable,   we can find in $A\setminus J$, a countable family of pairwise disjoint
subsets $V_{i}$ each of which is itself countably infinite. Index the members of
$V_i$ as $V_i = \brcs{v(i,1), v(i,2), \dots}$, where $v(i,n) \in
A\setminus J$. For each pair $(i,n)$, set $P_{i,n} = p_i +
2^{-n}v(i,n)/\| v(i,n)\|$. Then we note that $\lim_{n\to \infty} P_{i,n} =
p_i$, so that the countable set $P:= \brcs{P_{i,n}}_{\N \times \N}$ is dense.
 
Since the support of every element of $Q$ lies in $J_0$, $P \cap Q = \emptyset$;
since for any fixed pair $(i',n')$, the element $v(i',n')$ appears exactly once in the
supports of $P_{i,n}$ as $(i,n)$ varies, it easily follows that $P\cup Q$ is linearly independent (using linear independence of $Q$). 
\qed
 
We can ask whether 
\item{(*)}for every real infinite-dimensional Banach space, there exists
a dense linearly independent subset.
 
\noindent This has been solved affirmatively in [BDHMP; Proposition 3.2] (I am indebted to Ilijas Farrah for explaining their argument to me). 
 
The argument in the separable case used the fact that the dimension (as a
real vector space) exceeded the cardinality of one of the dense subsets.
However, this dimension property fails (when CH holds) for all
non-separable Banach spaces whose dimension is $\aleph_1$ (and probably
fails, in the presence of GCH, for all non-separable Banach spaces). The argument of the separable case can be adapted if $B$ has a dense subspace of codimension at least as large as its dimension. 
 
A weaker hypothesis that would still be enough for our purposes is the
following:
 
\item{(**)} for any real infinite dimensional Banach space, there exists a
linearly independent set $\brcs{e_{\alpha}}$ whose $\Z$-span (that is,
$\sum e_{\alpha}\Z$)
is dense.
 
\noindent Whenever this occurs, $\sum e_{\alpha}\Z$ is free  (on
$\brcs{e_{\alpha}}$) and has the property that every finitely generated
subgroup is discrete.
 
\Lem Corollary \thrtwo. Let $K$ be a metrizable Choquet simplex. Then there exists
a countable dense subgroup $G$ of $\Aff K$ that is free as an abelian group and for
which every finitely generated subgroup is discrete. Equipped with the
strict ordering, this $G$
is a countable simple dimension group with trace space $K$, with the
strong anti-finite dimensional property.
 
\Pf The Banach space $B = \Aff K$ is separable (since $K$ is metrizable),
and thus there exists a countable, dense, real linearly independent subset
of $B$, $\brcs{e_n}_{n=1}^{\infty}$; we  specify $Q= \brcs{ \pmb 1}$ (consisting of the constant function), so that if we take as order unit, $u = \pmb 1$, the affine representation agrees with the inclusion of $G$ in $\Aff K$. Then $G = \sum e_n \Z$ is dense; linear independence
over $\R$ implies that $G$ is free on $\brcs{e_n}$, and every finitely
generated subgroup of $G$ is contained in $\sum_{i=1}^m e_i \Z$, which is
discrete, and thus the subgroup is itself discrete. The rest is immediate.
\qed
 
These examples have a property that our original $R = \Z[x]$ does not: the group is free on a set which is itself dense.

By the affirmative solution to (*) in [BDHMP], we also obtain immediately the following. 

\Lem Proposition \thrthr. Let $K$ be a  Choquet simplex. Then there exists
a  dense subgroup $G$ of $\Aff K$ that is free as an abelian group and for
which every finitely generated subgroup is discrete. Equipped with the
strict ordering, this $G$
is a  simple dimension group with trace space $K$, with the
strong anti-finite dimensional property.

\SecT 4 Characterization of Anti-FD

The distinction between anti-finite dimensionality and its strong form is
not large.
 
\Lem Lemma \fouone. Let $(R,v)$ be an anti-finite dimensional approximately
divisible dimension group with order unit \st $\Inf R = \brcs{0}$. Let
$(G,u)$ be an approximately divisible dimension group with  order unit \st
$\partial_e S(G,u)$  is finite. Then there exists no continuous group
homomorphism $\Arrow \phi; G.R$  \st $\phi(G)$ contains an order unit of
$R$.
 
\Pf Suppose $\phi(G)$ contains an order unit. We may replace $G$ by its
simplification $G^s$; the topology is unchanged as is the set of
infinitesimals, but now $G^s$ is a simple dimension group, and thus so is
the quotient $G^s/\Inf (G^s)$. Since $\phi(\Inf G) \subset \Inf R $ (from
continuity of $\phi$ \wrt the pseudo-norms), $\phi$ induces a continuous
map from $G^s/\Inf G^s$ to $R$.
 
Hence we are reduced to the case that $G$ is a simple dimension group with
$\Inf G = 0$; the latter implies we may regard $G$ as a dense subgroup of
$\R^n = \Aff S(G,u)$. A continuous group homomorphism from $G$ is
automatically bounded, hence sends Cauchy sequences to Cauchy sequences,
so that $\phi$ extends to a map from the completion of $G$, $\R^n$, to the
completion of $R$, which is $\Aff S(R,v)$; call it $\Arrow \Phi; \R^n.
\Aff S(R,v)$.
 
Now $C:= \Phi^{-1}(\Aff S(R,v)^{++})$ is open; by hypothesis, it is
nonempty. Obviously, $C \cap -C = \emptyset$, $C-C = \R^n$ (since $C$, and
therefore $C-C$, contains an open ball), and $C + C \subseteq C$. Thus $C'
= C \cup \brcs{\pmb 0}$ is a proper convex cone in $\R^n$ containing an
open ball. Hence we may find a basis for $\R^n$, $\brcs{a_1, \dots, a_n}$
inside $C$. Let $D$ be the convex cone spanned by this basis; then $D$ is
a simplicial cone on $\R^n$ (that is, it is obtained from the original
ordering by applying an invertible matrix to everything in sight).

Density of  $G$  in $\R^n$ implies that if we impose the strict ordering on
$\R^n$ given by the basis (that is, if nonzero $x = \sum \lambda_i a_i$
where $\lambda_i \in \R$, then $x$ is in the positive cone iff $\lambda_i
> 0$ for all $i$), then $G$ becomes a simple dimension group \wrt this
orderingÑcall it (as an ordered group) $G'$. Since the positive cone is
contained in $D$, which in turn is contained in $C'$, it follows that
$\phi$ is positive as a group homomorphism from $G'$ to $R$. Since the
trace space of $G'$ is finite-dimensional, we reach a contradiction:
$\phi$ must be zero, as $R$ is anti-finite dimensional.
\qed

Now work toward completing the characterization of strong  anti-finite
dimensionality.
 
\Lem Lemma \foutwo. Let $\brcs{f_i}$ be a finite linearly independent subset of a
Banach space $B$. Then there exists $K > 0$ (depending only on
$\brcs{f_i}$) \st if $e_i$ are elements of $B$ with $\| e_i - f_i\| \leq
K$ for all $i$, then $\brcs{e_i}$ is linearly independent.
 
\Pf Form the finite-dimensional vector space $\sum f_i \R$; since all
norms on finite dimensional vector spaces are equivalent, choosing the
$l^1$-norm, there exists $K$ \st for all choices of $\lambda_i \in \R$, we
have $\| \sum \lambda_i f_i\| > K \sum |\lambda_i|$. If on the other hand,
$\sum \lambda_i e_i = 0$, then
$$
K \sum |\lambda_i| \geq \sum |\lambda_i| \| f_i - e\|  \geq \left\|\sum
\lambda_i(f_i - e_i) \right\| = \left\|\sum \lambda_i f_i \right\| > K \sum | \lambda_i |,
$$
a contradiction.
\qed
 
\Lem Lemma \fouthr. Suppose that $H$ is a finite rank subgroup of a Banach space
$B$. Then we may decompose $H = K \oplus F$ where $F$ is a discrete group,
$\overline H = \overline K \oplus F$ and $\overline K$
 is the maximal real subspace of $\overline H$.
 
\Pf Since $H\Q$ (the rational vector space spanned by $H$ inside $B$) is a
finite dimensional rational vector space, $\overline H$ is contained in a
finite dimensional subspace $V$ of $B$. Since all norms are equivalent on
finite dimensional spaces, we can write $\overline H =(\sum x_j \R) \oplus
(\sum f_i \Z)$
where the complete set $\brcs{x_j} \cup \brcs{f_i}$ is linearly
independent, and the maximal real subspace is the left summand.
Necessarily, $\sum f_i\Z$ is discrete.
 
By density, we may find $e_i$ in $H$ \st $\|e_i - f_i\|$ is as small as we
like; thus with $x_j = e_j$, we obtain that  $\brcs{x_j} \cup \brcs{e_i}$ is linearly independent; in particular,  $\sum e_i \Z$ is
discrete. It easily follows that $\overline H = (\sum
x_j \R ) \oplus (\sum e_i \Z)$.
 
Now consider the map $H \subset \overline H \to \sum e_i \Z$; this is onto
(since each $e_i \in H$ and the target is a free abelian group), and so
its kernel $K$ splits, that is, $H = K \oplus (\sum e_i \Z)$ where $K
\subseteq \sum x_j \R$. Now it is routine to verify that $\overline H =
\overline K\oplus (\sum e_i \Z)$, and thus $\overline K$ must be the
maximal real subspace of $\overline H$.
\qed
 
Lawrence [L] has shown that any countable discrete normed abelian group
(in particular, any countable discrete subgroup of a Banach space) is
free; this was extended by Steprans [S], dropping {\it countable.} In any
event, a finite rank discrete abelian group is free on a finite set, and this does not require either of these results.
 
\Lem Corollary \foufou. Let $(R,v)$ be an approximately divisible   dimension
group with order unit. 
\item{(a)} If $R$ is strongly anti-finite dimensional,
then $\Inf R = 0$ and  all finite rank subgroups of $R$ are discrete. In
particular, every countable subgroup of $R$ is free. 
\item{(b)}  If $R$ is 
anti-finite dimensional and $\Inf R = \brcs{0}$, then for every finite rank subgroup $H$ of $R$ whose closure is a vector space, $H \cap
R^{++} = \emptyset$.
 
\Pf First we prove that if $R$ is strongly anti-finite dimensional, then
$\Inf R = 0$. This is trivial: if $L =
\Inf R$, form the simple dimension group $H= \Q \oplus L$ with positive
cone $(q,l) > 0$ iff $q>0$. Then $H$ is a simple noncyclic dimension
group, and the map $(q,l) \mapsto l \in \Inf R$ is automatically
pseudo-norm continuous, hence must be zero. So $L$ is zero.
 
Let $H$ be a finite rank subgroup of $R$ that is not discrete. Taking the
closure in the Banach space $\Aff S(R,v)$, we have a split decomposition
$H = K \oplus F$ where $F$ is free and discrete and the closure of $K$ is
the maximal real subspace of $\overline H$; if $H$ is not discrete, then
$K \neq 0$. Since $K$ is dense in a finite dimensional vector space, we
can equip it with the structure of a noncyclic simple dimension group \st
the affine representation of $K$ is simply the embedding in $\overline K$.
Necessarily, the identity map $K \to K \subset R$ is continuous, so must
be zero, a contradiction to strong anti-finite dimensionality.
 
If instead, $R$ is only anti-finite dimensional, we can form $R^s$ and
factor out the infinitesimals; then the hypotheses still apply, so we can
reduce to the case that $\Inf R = 0$. Now using $K$ as in the preceding
paragraph, the hypotheses of lemma \fouone\ above yield a positive map from a
simple noncyclic dimension group to $R$, again a contradiction.
\qed

\Lem Corollary \foufiv. Let $(R,v)$ be an approximately divisible dimension group with order unit. 
Then $R$ is Anti-FD if and only if $\Inf R = 0$ and every finite rank subgroup of $R$ is discrete.

\Pf Follows from   \twoone\ and  \foufou. \qed
\comment
Necessity of the condition in (b) comes from the preceding. So assume the condition; there is no loss of generality in assuming  $\Inf R = 0$. Suppose $\Arrow \phi; G. R$ is a nonzero positive homomorphism from the approximately divisible dimension group with order unit $(G,u)$ \st $\partial_e S(G,u)$ is finite, to $R$. Since $R$ is simple, necessarily $v = \phi(u)$ is an order unit of $R$. If we replace $G$ by its simplification, $G^s$, then $\phi$ is still positive, and if factor out $\Inf G^s = \Inf G$, this persists. Hence we reduce to the case that $G$ is simple and $\Inf G = \brcs{0}$. In particular, we may regard $G$ as a subgroup of $\Aff S(G,u)$ which of course is a finite dimensional vector space, and the completion of $G$. 

Now $\phi$ extends to a continuous and positive map $\Arrow \Phi; \Aff S(G,u). \Aff S(R,v)$, and the image of $\Phi(G)$ is a finite dimensional subspace, in which $\phi(G)$ is dense. If $\phi(G)$ has finite rank, or more generally if $\phi(G)$ contains a finite rank subgroup which is dense in $\Phi(G)$, we are done. [In view of Example \twothr, dense subgroups of finite dimensional vector spaces need not have finite rank dense subgroups, so we are not finished.]

If $\phi(G)$ has infinite rank, we will obtain a contradiction. Since the closure of $\phi(G)$ is finite-dimensional, $\phi(G)$ spans (as a real vector space) its closure. Thus we may find a basis for the closure, $\brcs{p_1, \dots, p_n}$,  inside $\phi(G)$. Then $\phi(G) \subseteq L:= R \cap \sum p_i \R$
\endcomment

One would like a corresponding characterization of anti-fd, along the lines of, the simple dimension group $R$ is anti-fd if and only if every finite rank subgroup $H$ of $R/\Inf R$ whose closure is a vector space misses $R^{++}$. However, there is a problem arising from the phenomenon illustrated in Example \twothr, that dense subgroups of $\R^n$ need not contain any dense subgroups of finite rank. Necessary and sufficient, in case $R$ is simple,  is that (after factoring out $\Inf R$), if $H$ is a subgroup of $R$ whose closure is a finite dimensional real vector space, then $H \cap R^{++} = \brcs{0}$, but this is practically tautological. 

Finally we can give an example of a countable simple anti-fd with no infinitesimals that is not Anti-FD. Let $R = \Z[x] + (1-2x)\Q$ (a subgroup of $\Q[x]$) with the strict ordering from the interval $I= [\slfrac13,\slfrac23]$; this is a countable simple dimension group whose pure trace space can be identified with $I$. The set $\brcs{1, 1-2x,x^2, x^3, \dots}$ is linearly independent, and $(1-2x) \Q$ misses the positive cone. Let $\Arrow \phi; G. R$ be a positive group homomorphism, where $G$ has finitely many pure traces. 

If $\phi(G) \subset \Z[x] + n^{-1}(1-2x)\Z \subset (\slfrac 1n)\Z[x]$ for some $n$, then as the latter is order isomorphic to $\Z[x]$ with the strict ordering and is thus Anti-FD, and we obtain $\phi = 0$. 

Now $\phi(G)$ contains a real basis for its closure, say $\brcs{p_1,p_2, \dots, p_n}$, where each $p_i$ is a polynomial; this basis is contained in the finite dimensional space $V :=\sum_{j=0}^m x^i \R$ for some $m$, and $G$ is contained in this closed subspace. Hence the closure of $\phi(G)$ is contained in $V$. 

Since $\phi(G)$ thus consists of polynomials of degree less than or equal $m$ and is contained in $R$, that is, $\phi(G) \subseteq V \cap R$, and the latter is just  $\sum_{i=0,2,3,\dots, m} x^i \Z + (1-2x)\Q$Ñwhich is of finite rank. Hence $\phi(G) $ is of finite rank;  the closure of $\sum_{i=0,2,3,\dots, m} x^i \Z + (1-2x)\Q$ is just $(1-2x)\R +  \sum_{i=0,2,3,\dots, m} x^i \Z$; hence $\phi(G)$, being dense in a vector space, must map into the real subspace, so $\phi(G) \subseteq \overline{\phi(G)} \subseteq (1-2x)\R$. But this is impossible, as $(1-2x)\R$ contains no positive elements of $R$.

This example satisfies condition (b) of Proposition \twoone, but not (a). If instead we take $R = \Z[x] + \sqrt 2(1-2x)\Z$, the same argument shows that the resulting ordered group is again an anti-fd simple dimension group that is not Anti-FD, but satisfies condition (a) and  not (b). And $R = \Z[x] + (1-2x)\Q[\sqrt 2]$ is yet another example of a simple anti-fd but not Anti-FD dimension group, this time satisfying neither (a) nor (b).  

\SecT 5 Consequences

By direct translation from well-known results (too well-known to bother referring to), we have some consequences for unital AF C*-algebras and minimal actions on Cantor sets. Pro-fd simple dimension groups are $K_0$ of AF algebras which are direct limits of simple AF algebras each with finitely many pure traces. On the other hand, if $A$ is an AF algebra whose $K_0$ group is anti-fd, then $A$ contains no simple infinite dimensional AF-subalgebra with finitely many pure traces (more generally, the simple subalgebra need not be AF, but its K$_0$ group should be dense in a finite dimensional vector space).   

If $(X,T)$ is a minimal action on a Cantor set, and its dimension group (ordered \v Cech cohomology) is anti-fd, then $(X,T)$ is not even orbit equivalent to a minimal system with a non-trivial map to a minimal  action on a Cantor set that has only finitely many ergodic measures. 
 
\SecT 6 More examples

Very often, a vaguely ring-like structure on dimension groups is enough to guarantee the condition that every finite rank subgroup is discrete.

Let $A = \Z[x,x^{-1}]$ equipped with the coordinatewise ordering. For a $f \in A$, define $c(f)$, the {\it content of $f$}, as the greatest common divisor of the nonzero coefficients of $f$. Gauss' Lemma implies that $c(fg) = c(f) c(g) $ for $f,g \in A$. Let $(P_i)$ be a sequence of elements of $ A^+$, and form $\Cal R = \lim \Arrow \times P_i ; A. A$ as a partially ordered $A$-module. This is obviously a dimension group.

The order-theoretic properties of such dimension groups, beginning instead with $A = \R[x, x^{-1}]$ are studied in detail in [BH]. The determination of the positive cone and the pure trace (there called {\it state\/}) space is independent of the choice of coefficient rings ($\Z$ or $\R$), although of course the topological and underlying group structures are dependent on it. When we quote a result from [BH], it will apply to the pure trace space or to strong positivity. 

Define $Q_j = \prod_{i=1}^j P_i$. For $f \in A$, define $\Log f = \Set{k \in \Z}{(f,x^k ) \neq 0}$, where $(f,x^k)$ denotes the coefficient of $x^k$ in $f$.

Now let $R$ be the order ideal of $\Cal R$ generated by the constant function $1$; since an order ideal of a dimension group is again a dimension group, $R$ is a dimension group. We can identify $R$ with the following subgroup of $\Z[x^{-1} , p_i^{-1}]$,
$$
R \equiv R (p_i) = \Set{\frac f{Q_n}}{n \in \Z^+; \ \Log f \subseteq \Log Q_n}.
$$

Now assume that $c(p_i) =1$ for almost all $i$. We may delete finitely many $p_i$ and so assume that $c(p_i) = 1$ for all $i$. We now show that every finite rank subgroup of $R$ is discrete.

Form $G_n = \sum_{i \in \Log Q_n} \(\frac {x^i}{Q_n}\Z\)$; this is a finitely generated subgroup of $R$.

\noindent{(i)} {\it If $a \in R$ and $ma \in G_n$ for some nonzero integer $m$, then $a \in G_n$.} Write $a = f/Q_r$ for some $f \in A$ and $r \geq 1$ \st $\Log f \subseteq \Log Q_r$. Then there exists integers $t_i$ \st $m f/Q_r = \sum_{i\in \Log Q_n} t_i x^i/Q_n$. Let $g$ be the numerator of the right side, so that $\Log g \subseteq \Log Q_n$. Set $s = \max\brcs{r,n}$, and multiply the equation $m f/Q_r = g /Q_n$ by $Q_s$. Since $Q_s/Q_n$ and $Q_s/Q_r$ are products of $p_i$, each has content one. Thus $c(g) = c(mf) = m c(f)$. Thus $m$ divides $c(g)$, so that $h = g/m$ is in $A$. Thus $a = f/Q_r = h/Q_n \in G_n$.

\noindent{(ii)} {\it $G_n$ is discrete.} It suffices to show $\brcs{x^i/Q_n}{i \in \Log Q_n}$ is linearly independent over the reals. If $\brcs{\alpha_i}_{i \in \Log Q_n}$ are real numbers \st $\sum \alpha_i x^i/Q_n = 0$, on multiplying by $Q_n$, we obtain $\sum \alpha_i x^i = 0$, forcing all $\alpha_i = 0$.

\noindent {(iii)} {\it Any finite rank subgroup of $R$ is a subgroup of some $G_n$.}
If $S$ is a finite rank subgroup, we may find a finite subset $s_i $ of $S$ \st $S' = \sum (s_i \Z)$ has rank equalling that of $S$. There exists $n$ \st all $s_i \in G_n$. If $s \in S$, there exists a nonzero integer $m$ \st $m s \in S' \subset G_n$. By (i), $s \in G_n$. Hence $S \subseteq G_n$.

\noindent {(iv)} {\it Every finite rank subgroup  of $R$ is discrete.} Follows from (ii) and (iii).

Thus we have proved most of the following bifurcation result.

\Lem Proposition \sixone. Let $p_i$ be a family of elements of $A^+$, and set $R = R(p_i)$. Then all finite rank subgroups of $R$ are discrete if and only if $c(p_i) = 1$ for almost all $i$. If $c(p_i) >1$ for infinitely many $i$, then $ R \in \Cal H$.

\Pf If $c(p_i) = 1$ for almost all $i$, then we have just shown all finite rank subgroups are discrete. Otherwise, $c(p_i) > 1$ for infinitely many $i$. By telescoping, we may assume that $d_i:= c(p_i) >1$ for all $i$. Then we can write $p_i = d_i P_i$ for some $P_i \in A^+$, and it easily follows that if $U= \lim \Arrow \times d_i ; \Z.\Z$, then $R(p_i) \iso U \otimes R(P_i)$ as ordered abelian groups, hence we can write $R(p_i)$ as a limit of dimension groups of the form $U^{n(i)}$, hence $R(p_i) \in \Cal H$. 
\qed

We can also often tell when the ordered groups $R(p_i)$ have dense range in their affine representation (we normally use $1$ as the order unit for the representation). 

 Density of the image of $R(p_i)$ in its affine representation \wrt the order unit $1$ (hence as this is a dimension group, \wrt any order unit) is equivalent to the nonexistence of discrete pure traces. The kernel of a pure discrete trace of a dimension group is a maximal order ideal (with quotient $\Z$), by [GH, xxx]. Of course, we assume that all $p_i$ have at least two nonzero terms.
 
 There are two obvious maximal order ideals in $R(p_i)$, namely, the kernels of the point evaluations at $0$, that is, $\tau_0:f/Q_n \mapsto \lim_{t\downarrow 0}f(t)/Q_n(t)$ and $\tau_{\infty}:f/Q_n \mapsto \lim_{t\uparrow \infty}f(t)/Q_n(t)$. It is easy to evaluate the range of these two traces, in terms of leading and terminal coefficients. By multiplying each $p_i$ by a suitable power of $x$ (and converting the corresponding $P_n$), we may assume $0 = \min \Log p_i$ for all $i$, and let $d_i = \max \Log p_i$, so $D_n:= \max Log P_n = \sum_{i \leq n} d_i$. Then $\tau_0 (f/P_n) = (f,x^0)/(P_n,x^0)$ and $\tau_{\infty}(f/P_n) = (f,x^{D_n})/(P_n,x^{D(n)})$. 
 
 It easily follows that $\tau_0 (G) = \lim \Arrow \times (p_i, x^0); \Z.\Z$ and $\tau_{\infty}(G) =\lim \Arrow \times (p_i, x^{d^i}); \Z.\Z$. Hence necessary and sufficient for both these traces to have dense image is that 
 $$\eqalign{
 \text{for infinitely many $i$, }& (p_i, x^0) > 1\cr
 \text{for infinitely many $i$, }& (p_i, x^{d_i}) > 1.\cr
 }\tag\dag$$
 (Recall that we have altered the $p_i$ so the smallest exponent with nonzero coefficient is the constant term.)
 
 Hence, whenever $\ker \tau_0$ and $\ker \tau_{\infty}$ are the only maximal order ideals of $R(p_i)$, these conditions are also sufficient for density of the image of $R(p_i)$ in its affine representation, and thus for approximate divisibility.
 
 There are lots of situations in which these are the only maximal order ideals. We first make an observation: if $\Log p_i = \Log q_i$, then there is an obvious order isomorphism between the lattices of order ideals of $R(p_i)$ and those of $R(q_i)$. So to determine when they are the only maximal ideals, we can replace all the nonzero coefficients of all the $p_i$ by $1$. When we apply this process to $p_i$, the result will be denoted $p_i'$.
 Although the pure trace space changes when we go from $R(p_i)$ to $R(p_i')$, the lattice of order ideals does not.
 
 If $(p_i) $ is strongly positive (see [BH]); this is independent of the choice of coefficients], then the pure traces are precisely the point evaluations (including $\tau_0$ and $\tau_{\infty}$); the converse does not hold. Since maximal order ideals are automatically contained in kernels of pure traces (not necessarily discrete), it follows that if $(p_i)$ is strongly positive, then   $\ker \tau_0$ and $\ker \tau_{\infty}$ are the only maximal order ideals. Hence, if there exists a strongly positive sequence $(q_i)$ \st $\Log p_i = \Log q_i$, then $R(p_i)$ has only the two maximal order ideals. 
 
 If $p$ is a polynomial written in the form $p= \sum (p,x^i)x^i$, we call an exponent $i$ (or the corresponding monomial) an {\it isolani\/}%
 \plainfootnote{* }{based on Nimzovich's term for an isolated $d$-pawn in chess.}%
   if $(p,x^{i-1})$ and $(p,x^{i+1})$ are both zero. Next we note that if infinitely many $p_i$ have no leading or terminal isolani, then the sequence $(p_i')$ (obtained by replacing each nonzero coefficient in $p_i$ by $1$) is strongly positive (an easy consequence of the Superposition Lemma of [BH]). Hence if infinitely many $p_i$ have no isolani, then $R(p_i)$ has dense range in  its affine representation if and only if $(\dag)$ holds.
 
 The pure trace space is similarly just the set of point evaluations when infinitely many $p_i$ are  equal to each other (and projectively faithful, that is $\Log p_i - \Log p_i$ generates all of $\Z$). Hence if infinitely many $\Log p_i$ are equal to each other and projectively faithful, then $R(p_i)$ has only the two maximal order ideals, and so $\dag$ is necessary and sufficient for density of $R(p_i)$ in its affine representation.
 
 If there is a bound on the number of isolani, more generally, if there is an integer $N$ \st infinitely many $p_i$ have at most $N$ isolani and their diameter (distance from the smallest isolani to the largest) is bounded, then I think the pure trace space of $R(p_i')$ again is the set of point evaluations. The thesis by Alan Kelm [K] contains numerous results on strong positivity and related ideas, improving some of those in [BH].
 
 On the other hand, there exist  lacunary choices for $(p_i)$ for which there are lots of maximal order ideals, and even lots of discrete traces (the former depends only on $(\Log p_i)$, the latter depends on the actual coefficients). For example, let $p_i' = 1 + x^{3^i}$; these aren't projectively faithful, but we can also take $p_i' = 1 + x^{3^i} + x^{5^i}$. Obtaining conditions on the coefficients in lacunary examples so that density holds is not trivial. 
 
 Anyway, we have at least some results.

 \Lem Proposition \sixtwo. Suppose that $p_i \in A^+$ \st $|\Log p_i|\geq 2$ for all $i$ and  ($\dag$) holds. Then sufficient for $R(p_i)$ to have dense image in its affine representation is either
 \item{(a)} for infinitely many $i$, $p_i$ has no leading or terminal isolani.
 \item{(b)} there exists an infinite subset $S$ of $\N$  \st for all $i \in S$, $\Log p_i$ are equal to each other.
 
 If there is a bound on the degrees of the (nontrivial) polynomials $p_i$, then (b) applies. If there are no gaps, then there are no isolanis either (gaps can even persist, e.g., if $p_i' = (1+x)(1+ x^{i!})$, then $(p_i')$ is strongly positiveâ because of the factor $(1+x)$â but all finite products will have gaps),  and so the proposition applies to such sequences as well.

\Lem Corollary. \sixthr. Suppose that $p_i \in A^+$ \st $\Log p_i \geq 2$ for all $i$. If all of the following hold, then   $R(p_i)$ is Anti-FD.
\item{(i)} for infinitely many $n$,
 the coefficient of the smallest degree term in $p_n$ exceeds $1$;
\item{(ii)} for infinitely many $n$,
 the coefficient of the largest degree term in $p_n$ exceeds $1$;
\item{(iii)} for almost all $n$, $c(p_n) =1$
\item{(iv)} either (a) or (b) of Proposition \sixtwo.

\comment
In general, we cannot conclude that $R$ has dense image in its affine representation. The first, obvious obstruction is that if the leading or the terminating coefficients of all the $p_i$ are one, then $R$ admits a discrete trace (that is, a trace whose range is $\Z$). A second obstruction is that in this generality, the pure trace space is difficult to describe.

Fortunately, there are a lot of special cases, especially if $p_i$ are all equal, or if the $p_i$ have some sort of unimodality conditions (as in [BH]). First, in many cases, the pure trace space is the two-point compactification of $\R^{++}$, that is (up to homeomorphism, a closed interval). We note that in general, we have the following traces. If $r \in \R^{++}$, $\phi_r (f/Q_n) = f(r)/Q_n(r)$ (evaluation at $r$); if $r = 0$, $\phi_0 (f/Q_n) = \lim_{t\to 0} f(t)/Q_n(t)$ (that this exists uses l'H™pital's theorem and the fact that $\Log f \subseteq Q_n$); and if $r = \infty$, $\phi_{\infty} (f/Q_n) = \lim_{t\to \infty}$ (exists for the same reasons). These are traces, and in the weak topology, this set is the two-point compactification of the positive reals.

In general the point evaluation traces need not be pure, nor even if they are pure, they need not form a dense subset of the pure trace space. Fortunately, there are sufficient conditions for this set of pure traces to be exactly the pure trace space:

\item{(a)} If $\Log P_i = \Log P_j$ for all $j$ and $\brcs{(P_i, x^s)/(P_j,x^s)}_{i,j}$ is bounded for each $s \in \Log P_1$; this can be somewhat extended. A special case occurs when $P_i$ are all equal.

\item{(b)} When the sequence $\brcs{P_i}$ is strongly positive [BH]; that is, if $f \in \Z[x]$ and $f|\R^{++} > 0$, then for all $k$, there exist $m > n > k$ \st $P_m \cdot P_{m-1}\cdot \dots \cdot P_{n+1} \cdot P_n \cdot f$ has no negative coefficients.  Numerous sufficient conditions on the distributions of the coefficients of $P_i$ were given in[BH] and also in $[D]$ for this to hold (there, the $f$'s were allowed to be real polynomials, but this makes no difference). These mostly involved unimodality or strong unimodality of the coefficients of the polynomials, but there was also a superposition principle, which allowed non-unimodal sequences.

In addition, there exist non-strongly positive sequences for which the pure trace space is the two-point compactification of the positive reals [BH].

In the cases where we know that the pure trace space is $X = \R^{++} \cup \brcs{\pm\infty}$ (the $-\infty$ of course being the limit approaching zero), we can guarantee that $R$ has dense image in $C(X,\R)$, the affine representation. Since $R$ is a dimension group, sufficient for density is that $\tau(R)$ be noncyclic for every pure trace $\tau$ [GH, xxx].

Let $m(i) = \min\Set{j}{(p_i,x^j) \neq 0}$ and $M(i) = \max\Set{j}{(p_i,x^j) \neq 0}$. Then $\phi_0 (R) = \cup \(1/\prod_{i=1}^n(p_i,x^{m(i)})\Z\)$ and $\phi_{\infty}(R) = \cup \(1/\prod_{i=1}^n(p_i,x^{M(i)})\Z\)$. Hence $\phi_0(R)$ will be dense in $\R$ iff $(p_i, x^{m(i)})> 1$ for infinitely many $i$, and $\phi_{\infty}(R)$ will be dense in $\R$ iff $(p_i, x^{M(i)})> 1$ for infinitely many $i$. At the point evaluation traces, $\phi_r(R)$ is automatically dense if we assume simply that  $p_i$ is not a multiple of a monomial.

Hence if we know (a) that $X$ is the pure trace space of $R$, (b) that  $c(p_i) = 1$ for almost all $i$, and (c) that $(p_i,x^{m(i)}) > 1$ and $(p_i,x^{M(i)})> 1$ each for infinitely many $i$, then the image of $R$ is dense in its affine representation.

We have thus obtained the following.

\Lem Proposition. Suppose the $p_i \in \Z[x^{\pm1}]^+$ are nonconstant and satisfy $c(p_i) = 1$. Let $R$ be the order ideal generated by $1$ in $\lim \Arrow \times p_i; \Z[x^{\pm1}].\Z[x^{\pm 1}]$.
\item{(a)} Then every finite rank subgroup of $R$ is discrete.
\item{(b)} If both $(p_i, x^{m(i)}) > 1$ for infinitely many $i$ and $(p_i, x^{M(i)}) > 1$ for infinitely many $i$, {\it and\/} the pure trace space of $(R,1)$ is the two-point compactification of $\R^{++}$ (identified with point evaluations), then $R$ has dense range in its affine representation.{\par}
\noindent In particular, if $\brcs{p_i}$ is strongly positive, then the pure trace space of $(R,1)$ is the two-point compactification of $\R^{++}$. 
\endcomment
(Conditions (i) and (ii) together are equivalent to ($\dag$) for $R(p_i)$.)
Now we can take the simplifications of such $R$, so obtaining simple dimension groups that are Anti-FD.

For example, if $p_i = 2x+3$ for all $i$, all the conditions hold, and the corresponding $R$ (which also happens to be a ring, resembling the GICAR dimension group, but with values $\Z[1/3]$ and $\Z[1/2]$ instead of $\Z$ at the two endpoints. (Here $R $ is the ring generated by $Z= 1/(2x+3)$ with positive cone generated multiplicatively and additively by $Z$ and $xZ$. If $p_i = 3 + x + 2x^2$, then the corresponding ring is not singly generated, but all the results still apply.

We can actually extend this to several variables (all the relevant results proved in [H]), as well as ordered subrings arising from actions of compact groups (since discreteness of finite rank subgroups will be inherited, and in these cases, we know the pure trace space is typically the orbit space \wrt   the Weyl group of the pure trace space of the original).

Some of these are candidates to be initial objects for approximately divisible dimension groups.

\SecT 7 Initial object preliminaries

The proof of the following is based on the proof of [H, Proposition 1.2, last few paragraphs]. It was obtained {\it en passent\/} in the case of a partially ordered ring with $u =1$; however, the proof is identical, and does not require the ordered ring structure. The condition that $\widehat G$ contain $\hat u\R$ in its closure is equivalent to the stronger property of approximate divisibility when $G$ is a dimension group, but not for more general partially ordered groups.

\Lem Lemma \sevone. Suppose $(G,u)$ is a partially ordered unperforated abelian group with order unit that the closure of $\hat G$ contains $\hat u\cdot \R$. Let $p$ and $q$ be relatively prime positive integers; then   for all $\epsilon > 0$ for all $r \in (0,1]$, there exist order units $v$ and $w$ \st $u= pv + q w$ and  $\|\hat v - r\pmb 1/p\| < \epsilon/p$ and $\| \hat w - (1-r)\pmb 1/q\| < \epsilon/q$.

\Pf By interchanging $p$ and $q$ if necessary, we may assume that $r \in [\slfrac 12,1]$.  There exists a positive integer $p k \equiv 1 \mod q$, so that we can find integer $t$ \st $pk = qt + 1$. Now form the subgroup $ku\Z + qG$ of $G$; since this contains $qG$, its image in the affine representation (\wrt $u$) is dense. Hence there exists $z \in G$ \st 
$\pmb 1 \cdot \max\brcs{r-\epsilon, \slfrac 12}/p  < \widehat{(k u - q z)} < \pmb 1 \cdot r/p < 1/p$. Set $v = k u - q z$, so $\hat v \gg 0$, and thus $v$ is an order unit. Now $u - pv = u - pku + pqz$, so $u -pv = (1-pk)u + q z = q(t u + z) $. Set $w = tu +z$; then $u = pv + qw$; moreover, $\hat w = \widehat{u - pv}/q = q^{-1}(\pmb 1 - p\hat v) \gg 0$, the latter since $\hat v < p^{-1}\pmb 1$. Hence $w$ is an order unit.

It also follows that $\| \hat w - \pmb 1 \cdot (1-r)/q\| < \epsilon /q$ (just us the triangle inequality applied to $\| q \hat w - (1-r)\pmb \| = \|\pmb 1 - p\hat v - (1-r)\|$.
\qed

The following was proved by Perera and R\o rdam [PR], in the context of homomorphisms from finite dimensional algebras to real rank zero C*-algebras (it applies to AF algebras in particular). Here and in the preceding, we do not require $G$ to be a dimension group, merely to have dense image in its affine representation. If $G$ is a dimension group, then dense image is equivalent to all pure traces not being discrete.

\Lem Corollary \sevtwo. Let $(G,u)$ be a partially ordered abelian group with order unit \st the image of $G$ in its affine representation is dense. Suppose that $\brcs{p_i}$ is a finite set of positive integers \st $\gcd\brcs{p_i} = 1$. Given an order unit $U$ of $G$, there exist order units $v_i$ \st $U = \sum p_i v_i$.

\Pf Let $q = \gcd\brcs{p_2, \dots, p_n}$, so that $p_1$ and $q$ are relatively prime. Given an order unit $u$, by the preceding there exist order units $v_1$ and $w$ \st $u = p_1 v_1 + q w$. Now $qG$ has dense image in whatever affine representation $G$ has dense image. Since $\gcd\brcs{p_i/q}_{i=2}^n = 1$, so by induction applied to $qG$,  $qw$ can be expressed as as $\sum_{i\geq 2} p_i q^{-1} z_i$ where $z_i$ are order units in $qG$. Then $z_i = q v_i$ for some $v_i$ in $G$, necessarily order units, and $ u = \sum p_i v_i$, completing the induction.
\qed

In particular, if $\gcd\brcs{p_i} = 1$, then $\prod_{i=1}^n M_{p_i} \C$ is an initial object for approximately divisible AF algebras.
The paper by Perera and R\o rdam [PR] proves more; real rank zero and no representations which hits the compacts is sufficient to get this.

There is a peculiarity, noted previously, in the characterization of Anti-FD partially ordered groups. Another approach  emphasizes the peculiarity.

Let $(H,v)$ be an unperforated partially ordered abelian group with order unit $v$, and suppose that $\Inf H = \brcs{0}$. For $m$ a positive integer, we say $(H,v)$ (or just $H$) is {\it Anti-FD\paren{m}\/} if all subgroups $H_0$ of rank $m$ or less are discrete (\wrt the norm topology induced by the order unit $v$); that is, $\widehat{H_0}$ is a discrete subgroup of $\Aff S(H,v)$ in the sup-norm).

By the characterization above (Corollary \foufiv), $(H,v)$ is Anti-FD($m$) for all $m$ iff $(H,v)$ is Anti-FD. Examples include the critical groups of [H2]: a subgroup $G$ of $\R^m$ is {\it critical\/} if it is dense and of rank $k+1$; equipped with the strict ordering inherited from the vector space, these are simple dimension groups with $m$ pure traces, and every subgroup of rank $m$ or less is discrete. Hence Anti-FD$(m)$ simple dimension groups which are not Anti-FD$(m+1)$ exist for every $m \in \N$.

Anti-FD($m$) partially ordered groups have the expected (relatively weak) property dealing with continuous homomorphisms.

\Lem Proposition \sevthr. Suppose $(G,u)$ is an approximately divisible unperforated partially ordered abelian group, and $m$ is an integer \st $\rk G/\Inf G \leq m$. If $(H,v)$ is Anti-FD($m$) and $\Inf H = \brcs{0}$, then any continuous group homomorphism $\Arrow \phi; G.H$ is zero.

We require a few elementary lemmas. The norm (or pseudo-norm) on a group $(G,u)$ is determined by the order unit $u$ (changing the order unit changes the norm to an equivalent one, but even for simple dimension groups, may actually change the Choquet simplex into one that is not affinely homeomorphic to the original!).

\Lem Lemma \sevfou. Let $(G,u)$ be an unperforated approximately divisible group, and $n$ a positive integer. Let $S_n = \Set{g \in G^{++}}{\| \hat g \|  < \slfrac 1n}$. Then $S_n$ generates $G$ as an abelian group.

\Pf It suffices to show $G^{++}$ is contained in the group generated by $S_n$. Select an order unit $w$, and a relatively prime pair of  positive integers, $p$ and $q$. By Lemma \sevone\ above, there exist order units $y$ and $z$ in $G$ \st $w = p y + q z$. Since $w$, $y$, and $z$ are all in the positive cone, we have $\| \hat y \| \leq \| \hat w \|/p$ and $\|\hat z \| \leq \| w \| / q$. If we choose $p, q > n \|\hat w\|$, we have $y,z \in S_n$, and we are done.
\qed

\Lem Lemma \sevfiv. Let $(G,u)$ and $(H,v)$ be unperforated partially ordered groups with order unit, and let $\Arrow \phi;G.H$ be a positive group homomorphism (not necessarily sending
order units to order units). Then $|\widehat {\phi (g)}| \leq \| \hat g\| \cdot \| \widehat {\phi(u)}$ (that is, $\|\!| \phi \|\!| \leq \|\widehat{\phi(u)}\|$).

\Pf Suppose $g \in G$ and $\| \hat g \| = \delta > 0$. If $p$ and $q$ are positive integers \st $\delta < p/q$, then $\| \hat g \| < p/q$ entails $-q \pmb 1 \ll p\hat g \ll q \pmb 1$, so by unperforation, $-q u \leq p g \leq q u$. Applying $\phi$, we obtain
$- q \phi(u) \leq p \phi(g) \leq q \phi(u)$, and since $\phi(u)$ is positive, we deduce
$\| \widehat{\phi(g)}\| \leq (p/q)\| \widehat {\phi(u)}\|$. Since this is true for all $p/q > \delta$, the result follows.
\qed

\Pf (Proposition \sevthr). By definition, a continuous group homomorphism maps $\Inf G$ to $\Inf H = \brcs{0}$, so without loss of generality, we may factor out $\Inf G$, and thus assume that $G$ itself has rank $m$ or less. Then $H_0:= \phi(G)$ has rank at most $m$, and so is discrete in the affine representation of $(H,v)$. There thus exists $\epsilon > 0$ \st if  $\| \hat h \| < \epsilon$ for $h \in H_0$, then $\hat h = 0$, and thus $h = 0$.

From the continuity of $\phi$ at $0$, given $\epsilon$, there exists a positive integer $n$ \st for all $g \in G$, $\|\hat g\| < 1/n$ entails $\| \widehat{\phi(g)}\| < \epsilon$.
Thus if we define, as in \sevfou, $S_n = \brcs{g \in G}{\| \hat g \| < 1/n}$, then $\phi(S_n) = \brcs{0}$, and since $S_n$ generates $G$ as an abelian group, $\phi(G) = 0$.
\qed

Now we see the peculiarity. If we let $m\to \infty$, we only obtain that there are no nontrivial homomorphisms from finite rank simple dimension groups, whereas we know that we can replace finite rank by finite pure trace space, and we also know that not all finite pure trace space dimension groups can be written as a direct limit of finite rank approximately divisible groups. 

\SecT 8 Initial objects

Let $(H,v)$ be a partially ordered unperforated abelian group with order unit $u$. We say that $(H,v)$ is an {\it initial object\/} if for all approximately divisible partially ordered  abelian groups with order unit, $(G,u)$, there exists an order-preserving group homomorphism $(H,v) \to (G,u)$. Normally, $H$ will be a dimension group ($G$ need not be, but in cases of interest, it will be). This is designed to be compatible with initial objects in the class of unital $C*$-algebras (typically AF), via their ordered $K_0$-groups. Obviously if $v =2x$ for some $x \in H$ (and necessarily in $H^{++}$ as a consequence of unperforation), then $(H,v) $ cannot be an initial object, although it  is still possible for $(H,x)$ to be one. 

Initial objects in the category of dimension groups (which can be defined in many different, and equally plausible ways) that are approximately divisible  are automatically at least anti-fd (by [ER], the pure trace space of an approximately divisible initial object is infinite). Some of the $R(p_i)$ are initial objects; if $p_i = 1+x$ for all $i$, this was shown in [ER]. We show by completely different methods that if $p_i = a_i + b_i x$
with $a_i, b_i$ positive integers exceeding $1$ \st $\gcd(a_i,b_i) = 1$ and $\sum_i \min\brcs{a_i,b_i}/(a_i + b_i) < \infty$, then $R(p_i)$ is an initial object. Their pure trace spaces are all the one-point compactification of $\Z^+$.  All $R(p_i)$ are the ordered Grothendieck group of the fixed point algebras of product type actions of a circle on UHF C*-algebras. 

\comment
 Initial objects must have free K${}_0$? Initial objects are anti-fd, which means almost free. If it has infinitesimals, can it be initial object? Yes. Take a simple AF initial object, take its K$_0$ and  direct sum it with an arbitrary kernel (strict ordering), then there exists an AF with an algebra map, necessarily onto on K_0.
\endcomment

If $g, h \in G$ and $k$ is a positive integer, then we use either $ g\equiv h \mod k$\ \  or \ $ g\equiv h \mod kG$\ \   to denote $h -k \in kG$.

 \Lem Lemma \eigone. Let $(G,u)$ be a partially ordered abelian group which has dense image in $\Aff S(G,u)$. Then for all integers $k$, all $\epsilon > 0$, all $g \in G$ there exists $h \in G$ \st $ h \equiv g \mod k$ and $\|\hat h \| <\epsilon$.
 
 \Pf Since the image of $G$ is dense, so is that of $kG$ and therefore that of $-g + kG$. Hence $-g + kG$ contains elements of arbitrarily small norm.\qed

We will show the following. Let $\brcs{u_i}_{i=0}^{N-1}$ be a set of order units \st for some $\delta$, there exists a positive rational $c$ \st $\| \hat{u_i} - c\| < \delta$ for all $i$ (in this case, $c$ represents the constant function with value $c$). Suppose that $1 < a < b$ are positive and relatively prime integers. Then there exist elements $\brcs{v_i}_{i=0}^{N}$ of $G$ \st 
 $$\eqalign{
 \left\| \hat {v_i} - \frac c{a+b}\right\| & < \frac{\delta}{b-a}  \qquad \text{for $i = 0,1,2, \dots, N$}\cr
 u_i &= b v_i + a v_{i+1}  \qquad \text{for $i = 0,1,2, \dots, N-1$.}\cr
 }$$
 A consequence is that if $\delta < c(b-a)/(a+b)$, then the $v_i$ are all order units, hence in $G^+$. 
 
 We break this into three steps. We form $G_0 = G \otimes \Q$, elements of which we regard as formal fractions, $g/q$ where $g \in G$ and $q \in \Q^{++}$; it has the obvious positive cone. Since it is  a vector space over the rationals, solutions of linear systems (with integer coefficients) are relatively easy to deal with.
 
 We define the $N \times (N+1)$ matrix $A$ (with rows indexed $0,1,2,\dots, N-1$ and columns indexed $0,1,2,\dots, N$) via 
 $$
 A_{ij} = \cases b & \text{if $i=j$}\\
 a & \text{if $j = i+1$}\\
 0 & \text{else.}\\
 \endcases
 $$
 The equations in the second line above can be compressed into the simple $AV = U$, where $V = (v_i)^T$ and $U = (u_i)^T$. 
 
 For each $j = 0, 1,2, \dots, N$, define $T_j = u_0 - ab^{-1} u_1 + \dots + (ab^{-1})^j u_j$, that is, $T_j = \sum_{i=0}^j (\slfrac ab)^i u_i \in G_0$. 
 
 \Lem Lemma \eigtwo. All solutions $X = (x_i)^T \in G_0^{N+1}$ (indexed beginning at zero) to $AX = U$ are of the form given by 
 $$
 x_i = \(\frac {-b}a\)^i \(x_0 -  \frac {T_{i-1}}b\) \quad {i = 1,2,\dots, N}. 
 $$
 
 \Pf We first describe all the solutions to the homogeneous equation $AX = 0$. If $Z = (z_i)^T$ is a solution, then $b z_0 + a z_1	=0$ implies $z_1$ is a rational multiple of $z_0$, and by induction, every $z_i$ is a rational multiple of $z_0$. Hence we may write $Z = (1, q_1,q_2, \dots, q_n)^T z_0$ where all the $q$s are rational numbers. We easily see from the equations that $q_i = (-b/a)q_{i-1}$, and thus $z_i = (-b/a)^i z_0$. Let $w = (1, -b/a, b^2/a^2, \dots , (-b/a)^N )^T \in \Q^{N+1}$. We have just shown that all solutions to $AZ = 0$ are of the form $Z = w z_0$. 
 
 Next, define the vector $S = (S_i)^T$ where $S_i = (-b/a)^{i-1}T_{i-1}/a$ for $1 \leq i \leq N$ and $S_0= 0$. We claim that $AS = U$, so $S$ is a particular solution. The top coordinate of $AS$ is $aT_0/a = T_0 = u_0$. We simply calculate the rest; for $1 \leq j \leq N-1$,
 $$\eqalign{
 (AS)_j & = b S_{j} + a S_{j+1} \cr
 & = \(\frac {-b}{a}\)^{j-1}\( \frac ba T_{j-1} + \frac {-ba}{a^2} T_{j} \)\cr
 & =  \(\frac {-b}{a}\)^{j} \(T_j - T_{j-1}\) = \(\frac {-b}{a}\)^{j}\cdot \(\frac {-a}{b}\)^{j} u_j \cr
 & = u_j.\cr
 }$$
 
 Now let $X$ be any solution in $G_0^{N+1}$ to $AX = U$. Then $A(X-S) = \pmb 0$, so there exists $z_0 \in G_0$ \st $X = S + wz_0$, which is precisely the desired statement. (Note that $x_0 = z_0$.)\qed
 
 \Lem Lemma \eigthr. Suppose $u_i$ in $G_0$ satisfy $\| \hat u_i - c\| <  \delta$ for some positive constant $c$ and some $\delta > 0$. Then there exists a solution $X  =(x_i) \in G_0^{N+1}$ to $AX = U$ \st $\|\hat x_i - c/(a+b) \| < \delta/(b-a)$ for all $i = 0, 1,\dots, N$. 
 
 \Pf Set $Y = c\pmb 1/(a+b)$ (the column of constants). Then $AY = c\pmb 1$, which is close (in the infinity norm) to $U$. Let $M$ be the upper triangular $N \times N$ matrix (with coordinates indexed from $0$ to $N-1$) with $1$ directly above the diagonal and zeros elsewhere. Set $B = b\I + aM$ (so it constitutes the first $N$ columns of $A$). Then $B^{-1}$ (with entries in $\Q$) is $b^{-1}\(\I +  \sum_{i=1}^{m-1} (-a/b)^i M^i\)$, and we calculate $B^{-1}A$ has the identity as its first $N$ columns, and has $a((-a/b)^{N-1-i}))^T$ as the last column (recalling the indexing begins at $i=0$). 
 
 The equation $AX = U$ is equivalent to $B^{-1}AX = B^{-1}U$. Writing $X = xX_0^T$ (where $x$ is the final coordinate of $X$) and $B^{-1}U = (r_0,r_1, \dots, r_{N-1}^T$, we have $X_0 + ax ((-a/b)^{N-1-i}) = B^{-1}U$. We are free to vary $x$, and the value of $x$ uniquely determines $X$. Set $x = -cu/(a+b)$ (recall that the affine representation of $G$ is \wrt the order unit $u$, so $\hat x$ is  $ -c/(a+b)$, the constant function.

 The rest of the coordinates of $X = (x_i)$ ($x_N = -u/(a+b)$) are given by 
 $$
 x_i = \frac 1{b} \( u_{i} - \frac ab u_{i+1} + \(\frac ab\)^2 u_{i+2} + \dots + \(\frac ab\)^{N-1-i} u_{N-1} \) - \frac{acu}{a+b}\(\frac{-a}b\)^{N-1-i}.
 $$
 Now we approximate $\hat x_i$. 
 $$\eqalign{
 \left\| \hat x_i - \frac{c \(1 - (a/b) + (a/b)^2 + \dots + (-a/b)^{N-1-i} \)}b - \frac{ac}{a+b}
 \(\frac{-a}b\)^{N-1-i}\right\|  
 &\cr
< \frac {c\delta}b &\(1 + \frac ab + \frac {a^2}{b^2} + \dots + \(\frac{a}b\)^{N-1-i}\) \cr
 & =  \frac {c\delta}b \cdot \frac{1 - (a/b)^{N-i}}{b-a} b \cr
 & <  \frac {c\delta}{b-a}. \cr
 }$$
 
 Finally, $\frac{\(1 - (a/b) + (a/b)^2 + \dots + (-a/b)^{N-1-i} \)}b -\frac{a}{a+b}
 \(\frac{-a}b\)^{N-1-i} = 1/(a+b)$. 
 \qed
 
 \Lem Lemma \eigfou. If $X= (x_i)^T$ is a solution in $G_0^{N+1}$ to $AX= U$, then for all $\epsilon$, there exists a solution $V= (v_i)$ in $G^{N+1}$ to $AV = U$ \st $\| \hat x_i - \hat v_i\| < \epsilon $.
 
 \Pf First, we may find $y \in G$ \st $\|\hat x_0 - \hat y_0\| < \epsilon/a$. Now we perturb $y_0$ by an element of the form $h = \sum_{i=0}{N} h_i a^i$ where $\| \hat h_i \| < a^{-i -2}$ so that $v_i := (-b/a)^{i}(y + h - T_{i-1}/b) \in G$ (and $v_0 = y + h$). 
 
 Set $y_i = (-b/a)^i (y - T_{i-1}/b) \in G_0$. Multiply by $b^i$, so we obtain 
 $$
 a^i y_i = (-b)^i y  - (-b)^{i-1}T_{i-1}.
 $$
 The term $t_{i-1}:= (-b)^{i-1}T_{-1} = (-b)^{i-1}u_0 + (-b)^{i-2}a u_1 + (-b)^{i-3}a^2 u_2 + \dots + a^{i-1}u_{i-1}$ belongs to $G$. We want to find $h \in G$ with various properties \st $(-b)^i (y - h) \equiv - t_{i-1} \mod a^i$ for $i = 1, 2, \dots, N$. We write $h $ formally as $\sum_{i=0}^{N-1} a^i h_i$, and find the $h_i$ inductively. 
 
 If $i=1$, the equation reduces to $-b h_0 \equiv u_0 - b y \mod a$; since $\gcd\brcs{a,b} = 1$, this has a solution in $G$, and we may choose a solution $h_0 \in G$ \st  \st $\| h_0 \| < \epsilon/a^2$. Now suppose $h_0, \dots, h_{k-1}$ have been defined with $\| h_i \| < a^{-2i + 2}$ for $i = 0, 1,2, \dots, k-1$ and if $h^{(j)} = \sum_{i=0}^i a^i h_i$, then $(-b)^i (y - h^{(i-1)}) \equiv - t_{i-1} \mod a^i$ for all $ i \leq k$. We want to find $h_{k}$ \st $(-b)^{k+1}(y - h^{(k)} - a^k h_{k+1}) \equiv (-b)^k T_k \mod a^{k+1}$. 
 
 We have $T_k = T_{k-1} + (-a/b)^k u_k$, so $(-b)^k T_k = a^k u_k + (-b)(-b)^{k-1} T_{k-1}$. Also, the induction assumption implies $(-b)^{k}(y - h^{(k-1)}) \equiv (-b)(-b)^{k-1}T_{k-1} \mod a^k$; thus $x = ((-b)^{k}(y - h^{(k-1)})+ (b)(-b)^{k-1}T_{k-1}) /a^k$   belongs to $G$. The equation thus reduces to $a^k(-b)(x - (-b)^k h_k) \equiv 0 \mod a^{k+1} $. This is equivalent to $(-b) (x - (-b)^k h_k) \equiv 0 \mod 2$. Since $\gcd\brcs{a,b} = 1$, we can solve this with arbitrarily small norm. This completes the induction.
 \qed
 
 If we assume $b < a$ (and still have $\gcd\brcs{a,b} = 1$) instead of $a < b$, then we index the $u_i$ in reverse order; then the corresponding result (with $|b-a|$ replacing $b-a$) applies in this case.

 \Lem Theorem \eigfiv. Let $p_n = b_n + a_n x$ with $1 < a_n,  b_n$, $\sum \min\brcs{a_n,b_n}/(a_n + b_n) < \infty$, and $\gcd\brcs{a_n, b_n}=1$ for all $n = 1, 2, \dots $. Form $R(p_i) \equiv R = \Set{f/\prod_{i \leq n} p_i}{\Log f \subseteq \Log \prod_{i \leq n} p_i}$, the order ideal generated by $1$ in $\lim \Arrow \times p_n; A.A$. 
 Then $R$ is an initial object in the class of approximately divisible dimension groups. 
 
 \Pf Let $d = \prod_{i=1}^{\infty} |a_i - b_i|/(a_n + b_n)$; by hypothesis, $d > 0$. Let $(G,u)$ be an approximately divisible dimension group with order unit. Begin with $u^0 = u$. By \eigthr, we can find order units $u_0^1, u_1^1$ of $G$ \st $\left\| \hat u^1_i - 1/(a_1 +b_1) \right\| < d/|a_1 - b_1|$ and $u^0 = a_1 u_0^1 + b_1 u_1^1$. Since $d < |a_1 -b_1|/|a_1 + b_1|$, the $\hat u^1_i$ are strictly positive, and thus $u^1_i$ are order units, and the process may be continued. Now we proceed by induction. 

After the $N-1$st iteration, we can find $\brcs{u^N_j}_{j=0}^N$ \st $u^{N-1}_{i} = a_n u^{N}_i + b_n u^N_{i+1}$ for $i = 0,1,2, \dots, N-1$ and for all $j = 0,1,\dots, N$,
$$
\left\| \hat u^N_j - \frac 1{\prod_{t=1}^N (a_t + b_t)} \right\| < \frac d{\prod_{t=1}^N |a_t - b_t|} .
$$
Since $d < \prod_{t=1}^N |a_t - b_t|/(a_t + b_t)$, it follows that $\hat u_j^N$ are strictly positive, hence $u_j^N$ are order units, and we may continue to iterate this process. 

Now the assignment $x^j/\prod_{i \leq n} p_i \mapsto u^n_j$ (for $0 \leq j \leq n$) extends to a well-defined, positive homomorphism $R \to G$, sending $1$ to $u$. 
 \qed
 
 We can replace {\it approximately divisible dimension groups\/} by {\it partially ordered abelian groups with dense image in their affine function space,} since the constructions never use interpolation.
 
 The sequence $(p_n = b_n + a_n x)$ is either strongly positive (and more), in which case the pure trace space is the two-point compactification of $\R^{++}$, or not strongly positive, in which case the pure  trace space the one-point compactification of $\Z^+$. Strong positivity is equivalent (this is a very special case of a result in [BH]) in this case to $\sum \min \brcs{a_n, b_n}/(a_n + b_n) = \infty$, precisely the negation of the hypothesis in \eigfiv. When the sum converges (as in the last result),   the product $\prod p_i/(a_i + b_i)$ converges to an entire function, and the Maclaurin series coefficients determine the pure traces ([BH] again); the pure trace space is thus the one-point compactification of $\Z^+$, and by \sixtwo, the condition that $a_n , b_n > 1$ more than implies that $R$
 is approximately divisible.

\comment
 Suppose $p_n$ are positive polynomials (that is, their coefficients are nonnegative) and there are no gaps (that is, $(p,x^a)\cdot (p,x^{a+2}) \neq 0$ implies $(p,x^{a+1}) \neq 0$. Then the order ideal structure of $R$ (order ideal generated by $1$) is natrually isomorphic to that obtained from setting $p_n' $ to be $p_n$ with all nonzero coefficients set equal to one. This is trivial. Hence the only maximal order ideals are exactly the kernels of the evaluations at zero and infinity. Suppose an $R$ of this form admitted a pure trace $t$ \st $t(R)$ is discrete (such a $t$ is called a {\it discrete trace\/} in [G]]. By [GH, xxx], $\ker t$ is a maximal order ideal. Hence $\ker t$ would have to be the kernel of  evaluation at 0 or $\infty$, and thus $t$ itself would have to be evaluation at zero or infinity.
 
 Hence if we exclude evaluation at zero or infinity from being discrete traces, every pure trace of $R$ will have dense range in the reals, so by [GH, xxx], $R$ will have dense range in its affine representation. Obviously necessary and sufficient for the evaluations at 0 and $\infty$ not to be discrete are:{\par}
 \item{(*)} $(p_n, x^{m(n)}) > 1$ for infinitely many $n$
 \item{(**)} $(p_n, x^{M(n)}) > 1$ for infinitely many $m$
 where $m(n)$ is the least exponent $a$ \st $(p_n, x^a) \neq 0$ and $M(n)$ is the greatest exponent $a$ \st $(p_n, x^a) > 0$. 
 
 In particular, if (*) and (**) hold for our $p_n = a_n x + b$ (and $\gcd(a_n,b_n) = 1$, then the corresponding $R$ has dense image in its affine representation, even though the pure trace space is the two-point compactification of $\N$ (or one-point compactification of $\Z^+$ if we include evaluation at zero—the normalized constant coefficient—among the Maclaurin-type traces arising from convergence of the normalized infinite product, see [BH, last section]). 
\endcomment
 
 Here is an amusing class of initial objects, easily proved to be so, and easily determined when they have dense image in their affine representation. Moreover, they also provide examples of $R(p_i)$ for which checking merely at the two end evaluations is not sufficient to determine density; in fact, the collection of maximal order ideals is indexed by a Cantor space (in a natural way; in fact, the indexing is topological, if we use the usual topology on maximal order ideals), and we have to check density of the image at every one of the extremal traces.

 Let $p_i$ be a family of polynomials in $A^+$ \st $|\Log p_i| \geq 2$ for all $i$. We say the sequence (or the corresponding dimension group) $(p_i)$ is {\it non-interactive\/} if whenever $x^m \in \prod p_i$ (for some finite product of distinct $p_i$), there exist unique $w_i \in \Log p_i$ \st $m = \sum w_i$. The simplest example arises when $p_i = a_i + b_i x^{2^i}$ for $a_i, b_i \in \N$, but there are lots of others. The term comes from ergodic theory, actions of $\Z$, analysed \wrt measure-theoretic equivalence. 
 
 \Lem Proposition \eigsix. Let $(p_i)$ be a non-interactive sequence of polynomials in $A^+$ \st $c(p_i) = 1$ for almost all $i$. Then $R(p_i)$ is an initial object in the category of approximately divisible unital partially ordered unperforated abelian groups. The range of every pure trace of $R(p_i)$ is a subgroup of the rationals. Moreover, $(R(p_i),1)$ has dense image in its affine representation if and only if for infinitely many $i$, none of the coefficients  appearing $p_i$ is $1$ (that is, $(p_i, x^j) \neq 1$ for all $j$).
 
 \Pf Define $P_n =  \prod_{i \leq n} p_i$. Non-interactive means the Bratteli diagram for $R(p_i)$ (and also for $S(p_i)$, from which the former can be obtained by beginning at a single point in the top level) is  non-interactive, meaning if $j \neq k \in \Log P_n$, then $j + \Log p_{n+1} \cap k + \Log p_{n+1} = \emptyset$. Since the set of entries of $\Log p_{n+1}$ has greatest common divisor $1$, for each vertex of the diagram, corresponding to $x^i/P_n$ with $i \in \Log P_n$, we can solve the equations $U = \sum (p_{n+1},x^j)u_j$ in an arbitrary approximately divisible etc group $G$ (no matter what the order unit $U$ is) at each index, independently of the other indices on the same level. 
 
 If we use the notation $u_{i,n}$ to denote the order unit corresponding to $i \in \Log P_n$, we have found order units $v_{i + t,n+1}$ order units with $t \in \Log p_{n+1}$
 \st $u_{i,n} = \sum (p_{n+1},x^t)v_{i+t, n+1}$. Then for $s \in \Log P_{n+1}$ we define $(u_{s,n+1})$ to be $v_{i+(s-i),n+1}$ where $i$ is the unique element of $\Log P_n$ 
 \st $s \in i + \Log p_{n+1}$ (so necessarily $s-i \in \Log p_{n+1}$). Then the map $x^i/P_n \mapsto v_{i,n}$ (for $i \in \Log P_n$) yields a well-defined, positive group homomorphism $(R(p_i),1) \mapsto (G,u_0)$. Hence $R(p_i)$ is an initial object.
 
 Next, we show that the advertised condition is precisely what is needed to guarantee that $R(p_i) $ has no (pure) discrete traces. It is practically tautological that the the pure trace space can be identified with the path space of the reduced Bratteli diagram (obtained by collapsing all  multiple edges to single edges), which is a Cantor set. Along each path, the corresponding trace divides the value by the coefficient. Hence if all the $p_i$ admitted $1$ as a coefficient, there would be a path corresponding to repeated division by $1$, hence the image of the trace would be $\Z$. Conversely, if every path ran into infinitely many levels where all the nonzero coefficients of $p_n$ are not one, we would obtain division by infinitely many integers exceeding $1$, so that the range of the trace would be a noncyclic subgroup of the rationals. In that case, no extreme trace is discrete, hence $R(p_i)$ has dense range in its affine representation. 
 \qed
 
 So if we set $p_i = a_i + b_i x^{2^i}$, then  $(p_i)$ is non-interactive, and $R(p_i)$ is an initial object. It will have dense range in its affine representation if and only if  for infinitely many $i$, both $ a_i , b_i > 1$ hold. This is a stronger condition than the condition considered previously, which amounts to (in this case), $a_i > 1 $ for infinitely many $i$ and $b_j > 1$ for infinitely many $j$. These two conditions merely guarantee that the $\tau_0$ and $\tau_{\infty}$ are not discrete. These traces correspond to the extreme left path and the extreme right path in the Bratteli diagram, respectively. The case that $a_i = 2$ and $b_i = 3$ yields the ordered K$_0$-group of the AF-algebra given as the infinite tensor product,  $\otimes (M_2 \oplus M_3)$, considered in [ER]. 
 
There is an easy generalization of a sequence of non-interactive
polynomials, which also yields a family of initial objects.
 
Let $X$ be an infinite tree, and let $X^0$ be the Bratteli diagram
obtained from $X$ by labelling the edges with positive integers. We
sometimes identify $X$ with its path space. Let $G$ be the resulting
direct limit. The pure traces are given by paths in $X$ (not $X^0$) as
follows. Let $X_n$ be the set of vertices at the $n$th level, so $X_0 =
\brcs{x_0}$ consists of the initial point, and $X = \cup X_n$. Let $p:=
(x_0, x_1, \dots)$ be a path in $X$ (where $x_i \in X_i$), and suppose the
multiplicity of the edge $x_i \to x_{i+1}$ is $m(i)$. Let $\Arrow f; X_n .
\Z$; the trace $\tau_p$ is given by the map $[f,n] \to
f(x_n)/m(0)m(1)\cdot \dots \cdot m(n-1)$. It is easy to check that this is
well-defined and a pure trace.
 
Associated to each vertex  $x \in X$ is its vector of multiplicities,
$v(x)$; this  is the list of multiplicities on the edges emanating from $x$.
 
Moreover, $\ker \tau_p$ is a maximal order ideal (the order ideals of
$C(X,\Z)$ are in bijection with the order ideals of $G$), and thus
$\tau_p$ is pure. Moreover, $\brcs{\tau_p}_{p \in X}$ is a compact set of
pure traces (routine) homeomorphic to the path space, $X$. To check that
this is the pure trace space of $G$, we note that the embedding $\Z^{X_n}
\to \Z^{X_{n+1}}$ (of which $G$ is the direct limit) is an order
embedding:  $[f,n] \geq 0$ iff $f \geq 0$. But if $f$ is nonnegative along
every path, then its values at each of the points in $X_n$ are
nonnegative. Hence   $\brcs{\tau_p}_{p \in X}$ is a compact set of pure
traces that determines the ordering, and thus its closed convex hull is
the normalized trace space of $X$. It follows that $\brcs{\tau_p}_{p \in
X}$ is the pure trace space of $G$.
 
Now $G$ is an archimedean (in the strong sense) subgroup of $C(X,\R)$, and
density is equivalent to every $\tau_p$ having dense range. Sufficient for
this is that there exist infinitely many $n$ \st for every $x \in X_n$, $1
\not\in v(x)$.
 
It is a direct consequence of \sixtwo\ that $G$ is an initial object if $c(v(x) )= 1$ for all $x \in X$. 
Thus we have the following.

\Lem Proposition \eigsev. Let $X$ be an infinite  tree with root $x_0$, and $X_0$ the corresponding Bratteli diagram obtained by attaching positive-integer-valued weights to each vertex, and let $(H, u:=\chi_{\brcs{x_0}})$ be the resulting dimension group with order unit. 
\item{(i)} Sufficient for $(H,u)$ to be an initial object is that for all $x \in X$, $c(v(x)) =1$. 
\item{(ii)}Sufficient for $H$ to be approximately divisible is that for infinitely many $n$, for all $x \in X_n$, no entry of $v(x)$ is $1$. 
\item{(iii)} The pure trace space of $(H,u)$ is naturally homeomorphic to the path space of $X$. 
 
A reasonable question is whether all Anti-FD dimension groups are initial objects. If we extend the definition of initial object to require a one to one map  (as was proved for the Pascal's triangle example in [ER]), then such initial objects must be countable anf free (as abelian groups). So this stronger property excludes the non-Anti-FD but anti-fd example $\Z[x] + (2x-1)\Q$ discussed at the end of section four.

\long\def\Rf[#1] #2, #3. #4\par%
{\vskip 2pt \itemitem{[#1]} #2, {\it #3,} #4\par\vskip2pt}

\SecT References

\Rf [BDHMP] {T Bartoszynski, M Dzamonja, L Halbeisen, E Martinov‡, A Plichko},
On bases in Banach spaces. Studia Mathematica, 170  (2005) 147--171.

\Rf [BH] BM Baker \& DE Handelman, Positive polynomials and time dependent integer-valued random variables.
Canadian J Math 44 (1992) 3Ð41.

\Rf [BeH] S Bezuglyi \& D Handelman, Measures on Cantor set{\/\rm:} the good, the ugly, the bad.
Trans Amer Math Soc (to appear).

\Rf [C] I Chlodovsky, Une rŽmarque sur la reprŽsentation des fonctions continues par des polyn™mes ˆ coefficients entiers. Mat Sb 32 (1925) 472--475.

\Rf [EHS] {EG Effros, David Handelman, \& Chao-Liang Shen}, Dimension groups and their affine representations. Amer J Math 102 (1980) 385--407. 

\Rf [ER] GA Elliott \& M R\o rdam, Perturbation of Hausdorff moment sequences and an application to the theory of C*-Algebras of real rank zero. Operator Algebras
Abel Symposia, Volume 1 (2006) 97--115. 

\Rf [F] Le Baron O Ferguson, Approximation by polynomials with integral coefficients. Amer Math Soc, Rhode Island, 1980.

\Rf [G] KR Goodearl, Partially ordered abelian groups with interpolation. Mathematical Surveys and Monographs, 20, American Mathematical Society, Providence RI, 1986.

\Rf [GH] KR Goodearl \& David Handelman, Metric completions of partially ordered abelian groups. Indiana Univ J Math 29 (1980) 861--895.

\Rf [Gr] P  Griffith, Infinite abelian group theory. University of Chicago press, 1970.

\Rf [H] D Handelman, Iterated multiplication of characters of compact connected Lie groups. J of Algebra 173 (1995) 67--96.

\Rf [H2] D Handelman, Free rank $n+1$ dense subgroups of $\text{\/\bf R}^{n}$ and their endomorphisms.  J Funct Anal  46  (1982), no\. 1, 1--27.

\Rf [K] A Kelm, Strong positivity results for polynomials of bounded degree.  PhD thesis, University of Ottawa (1993).

\Rf [L] J Lawrence, Countable abelian groups with a discrete norm are free. Proc Amer Math Soc 90 (1984) 352--354.

\Rf [LL] AJ Lazar and J Lindenstrauss, Banach spaces whose duals are L${}^1$ spaces and their representing matrices. Acta Math 126 (1971) 165Ð193. 

\Rf [PR] F Perera and M R\o rdam, AF-embeddings into C*-algebras of real rank zero. J Funct Anal 217 (2004)  142--170.

\Rf [S] J Steprans, A characterization of free abelian groups. Proc Amer Math Soc 93 (1985) 347--349.

\vskip 10pt

\noindent Mathematics Dept, University of Ottawa, Ottawa K1N 6N5 ON, Canada; dehsg\@uottawa.ca

\comment
\Lem Example.  (continued) If $G$ is an approximately divisible dimension group with finitely many pure traces, there is no positive group homomorphism $\Arrow \alpha; G.R$.

\noindent Let $\Arrow \alpha; G.R$ be such a map. Without loss of generality, we may assume $G$ is simple, by imposing the strict ordering from its affine represention on it (this is coarser than the original, so $\alpha$ is still positive). We may also factor out the infinitesimals from $G$, since $R$ has none, and thus $\alpha$ would kill them anyway. In  the quotient ordering, the result is still a simple dimension group, and now the affine representation is just an embedding of $G$ (the quotient, the original $G$ modulo its infinitesimal subgroup) as a dense subgroup of $\R^n$. So we are reduced to the case that
$G$ is an dense subgroup of $\R^n$, with the strict ordering therefrom, and moreover, this is its affine representation.

Now $\alpha (G)$ is a subgroup of the free group $\Z[x]$, so is of course free. If $\alpha(G)$ were of finite rank, we could use the earlier argument in the finite rank case, that is, it would be finitely generated, and thus contained in a discrete group $\sum_{i=0}^n x^i \Z$ (for some $n$), and the bounded linear functional argument (which would work without reference to functionals, using bounded group homomorphisms instead) succeeds. Assume $\alpha(G)$ is of infinite rank.

Any discrete subgroup of $G$ has rank at most $n$ (since a discrete subgroup of $\R^n$ is free on an $\R$-linearly independent set). Also, any finitely generated subgroup of $G$ is discrete in $C([1/2,3/4],\R)$, as it is contained in $\sum_{i=0}^n x^i\Z$ for some $n$,  and $\brcs{1,x,x^2, \dots, x^n}$ is a real linearly independent subset of   $C([1/2,3/4], \R)$ (a polynomial that vanishes on an interval is zero). Since $\alpha(G)$ is assumed of infinite rank, it contains a subgroup $J$ of rank  $n+1$ (always true; even easier here, since $\alpha(G)$ is free, and thus has a $\Z$-basis). Being free and discrete, we may write it in the form $J = \oplus_{k=1}^{n+1} j_k \Z$. Thus there exist bounded linear functionals $s_i$ on $C([1/2,3/4], \R)$, $i=1,\dots, n+1$ \st $s_i (j_k ) = \delta_{ik}$ (the Kronecker delta). Each $s_i = \gamma_i - \beta_i$ where $\beta_i$ and $\gamma_i$ are (unnormalized) positive linear functionals on $C[1/2,3/4]$, hence traces on $\Z[x]$.

We define $B_i = \beta_i \circ \alpha$ and $C_i = \gamma_i \circ \gamma$, so that each is a trace on $G$, and thus $\sigma_i = B_i-C_i$ is a continuous group homomorphism (i.e., extends continuously to the affine representation of $G$, which is simply the completion, $\R^n$). However, if $g_j \in G$ are such that $\alpha(g_j) = k_j$, then we have that
$\brcs{\sigma_i}$ is linearly independent (since the set is dual to $\brcs{g_j}$). This means the trace space of $G$ has affine dimension at least $n+1$, a contradiction. \qed

A particular outcome is that $R$, a simple dimension group with pure trace space an interval, cannot be realized as a direct limit of approximately finite dimension groups each with only finitely many pure traces.
\endcomment

\end